\documentclass[12pt,oneside]{amsart}
\usepackage{amssymb,amscd,euscript}
\input diagrams.sty

\makeatletter \@addtoreset{equation}{section}
\@addtoreset{footnote}{section}\makeatother

\newcommand{\xref}[1]{{\rm \ref{#1}}}

\newcommand{\lins}[1]{\left|{#1}\right|}
\newcommand{\eqdef}{\mathrel{\,\stackrel{\mathrm{def}}{=}\,}}

\newcommand{\ov}[1]{\overline{#1}}
\newcommand{\comp}{\mathrel{\,\scriptstyle{\circ}\,}}

\newcommand{\1}{^{-1}}

\newcommand{\bDiv}{{\operatorname{BDiv}}}
\newcommand{\bCDiv}{{\operatorname{BCDiv}}}
\newcommand{\bir}{\dashrightarrow}
\newcommand{\ep}{\varepsilon}
\newcommand{\var}{\varphi}
\newcommand{\mt}[1]{\operatorname{#1}}
\newcommand{\Supp}{\operatorname{Supp}}
\newcommand{\Center}{\operatorname{Center}}
\newcommand{\Max}{\operatorname{Max}}
\newcommand{\Min}{\operatorname{Min}}
\newcommand{\Fix}{\operatorname{Fix}}
\newcommand{\Pic}{\operatorname{Pic}}
\newcommand{\Div}{\operatorname{Div}}
\newcommand{\Bs}{\operatorname{Bs}}
\newcommand{\Mov}{\operatorname{Mov}}
\newcommand{\Mv}{\operatorname{Mv}}

\newcommand{\N}{\operatorname{N}}
\newcommand{\NS}{\operatorname{NS}}

\newcommand{\mult}{\operatorname{mult}}
\newcommand{\ord}{\operatorname{ord}}
\newcommand{\sect}{\operatorname{sec}}
\newcommand{\Const}{\operatorname{Const}}
\newcommand{\CKM}{\operatorname{\scriptscriptstyle{CKM}}}
\newcommand{\s}{\operatorname{s}}
\newcommand{\m}{\operatorname{m}}
\newcommand{\e}{\operatorname{e}}

\newcommand{\PP}{{\mathbb P}}
\newcommand{\CC}{{\mathbb C}}
\newcommand{\QQ}{{\mathbb Q}}
\newcommand{\NN}{{\mathbb N}}
\newcommand{\RR}{{\mathbb R}}
\newcommand{\ZZ}{{\mathbb Z}}
\newcommand{\FF}{{\mathbb F}}

\newcommand{\F}{{\mt{f}}}
\newcommand{\gF}{{\mt{gf}}}

\newcommand{\sD}{\mathcal D}
\newcommand{\sE}{\mathcal E}

\newcommand{\sH}{\mathcal H}

\newcommand{\sP}{\EuScript P}

\newcommand{\sL}{\mathcal L}
\newcommand{\sM}{\mathcal M}
\newcommand{\sN}{\EuScript N}
\newcommand{\Oh}{\mathcal O}
\newcommand{\OOO}{\mathcal O}
\newcommand{\sR}{\mathcal R}
\newcommand{\sS}{\mathcal S}

\newcommand{\sG}{\mathcal G}

\renewcommand{\emptyset}{\varnothing}

\newcommand{\qq}{\mathbin{\sim_{\scriptscriptstyle{\QQ}}}}
\newcommand{\qr}{\mathbin{\sim_{\scriptscriptstyle{\RR}}}}

\newcommand{\down}[1]{\left\lfloor #1\right\rfloor}

\newtheorem{theorem}[equation]{Theorem}
\newtheorem{proposition}[equation]{Proposition}
\newtheorem{lemma}[equation]{Lemma}
\newtheorem{corollary}[equation]{Corollary}
\newtheorem{claim}[equation]{Claim}

\theoremstyle{definition}
\newtheorem{definition}[equation]{Definition}
\newtheorem{example}[equation]{Example}
\newtheorem{examples}[equation]{Examples}
\newtheorem{pusto}[equation]{}

\theoremstyle{remark}
\newtheorem{remark}[equation]{Remark}
\newtheorem*{remark*}{Remark}
\newtheorem{remarks}[equation]{Remarks}

\newenvironment{outline}{\begin{proof}[Sketch of the proof]}{\end{proof}}

\newcommand{\PLF}[1]{{\cite[#1]{PLF}}}
\newcommand{\PL}{{\cite{PLF}}}
\newcommand{\Isk}{{\cite{I}}}

\title{On Zariski decomposition problem}

\author{Yuri G. Prokhorov}

\thanks
{This work was is carried out  under the support of  grants
Leading  Scientific Schools 00-15-96085,  RFFI 02-01-00441, and
INTAS-OPEN 2000-269}

\address{
Department of Algebra, Faculty of Mathematics, Moscow State
Lomonosov University, Leninskie Gory, 117234 Moscow, Russia}

\email{prokhoro@mech.math.msu.su}

\date{}

\begin{document}

\begin{abstract}
We discuss different generalizations of Zariski decomposition,
relations between them and connections with finite generation of
divisorial algebras.
\end{abstract}

\maketitle\tableofcontents

The purpose of this survey is to clarify some details in \S\S 3-4
of \PL\ related to connections between different generalizations
of Zariski decomposition and finite generation of divisorial
algebras. The survey almost does not contain new results. All the
discussed results were known earlier (see \cite{Be}, \cite{Cu},
\cite{F3}, \cite{Ka}, \cite{Mo1}, \cite{Mo2}, \cite{Nak}, \PL,
\cite{Ts1}, \cite{Ts2}).

\section{Some preliminary facts}
All varieties are assumed to be normal, projective and defined
over the field of complex numbers $\CC$. However, almost all
results can be generalized to the relative situation $X/Z$.
Everywhere, if is not is indicated converse by a divisor we mean
an $\RR$-divisor, i.e., an $\RR$-linear combination of prime Weil
divisors. For any divisor $D$, we put
\[
\Oh_X (D)=\Oh_X (\down D),\qquad\lins{D}=\{F\mid F\sim D,\quad
f\ge 0\}.
\]
For definitions and main properties of b-divisors we refer to
Iskovskikh paper \Isk.

Recall that a divisor $D$ is said to be \emph{big,} if $\kappa (X,
D)=\dim X$, where
\[
\kappa (X, D)\eqdef
\begin{cases}
-\infty, &\text{if $\lins{nD}=\emptyset$ for all $n\in\NN$},
\\
\max\limits_{n\in\NN}\dim\Phi_{\lins{nD}} (X), &\text{otherwise},
\end{cases}
\]
is the \emph{Iitaka $D$-dimension} of $(X, D)$. A divisor $D$ is
said to be (\emph{semi}) \emph{ample} if $D=\sum\alpha_i H_i$,
where $H_i$ are integral (semi) ample divisors and
$\alpha_i\in\RR_{> 0}$. By the Kleiman criterion, the ampleness is
equivalent to the positiveness of the divisor on $\ov{NE}
(X)\setminus\{0\}$. The semiampleness is equivalent to the
existence of a contraction $\alpha\colon X\to Y$ and an ample
divisor $H$ on $Y$ such that $D\qr\alpha^*H$.

Recall that a divisor $F$ on a variety $X$ is said to be
\emph{b-semiample} (\emph{b-nef}) if there exists a model $Y$
dominating $X$ and an $\RR$-Cartier divisor $L$ on $Y$ such that
the b-divisor $\ov L$ is b-semiample (respectively, b-nef) and $F=
(\ov L)_X$. Equivalent: there exists a birational contraction
$f\colon Y\to X$ and a semiample (respectively, nef) $\RR$-Cartier
divisor $L$ on $Y$ such that $f_*L=F$. On a nonsingular variety
the b-semiampleness and semiampleness (respectively, b-nef and nef
properties) are equivalent \cite[Th. 6.1] {Z}.

\begin{pusto}
A \emph{rational} $1$-\emph{contraction} is a dominant rational
map $\alpha\colon X\dasharrow Y$ with connected fibers such that
\[
\dim\Center_Y G <\dim G\quad (=\dim X-1)
\]
for any prime exceptional b-divisor $G$ on $X$.
\end{pusto}

According to Hironaka for a rational $1$-contraction $\alpha\colon
X\dasharrow Y$ there exists a hut
\begin{equation}
\label{diag-bss}
\begin{diagram}
&&W&&
\\
&\ldTo^{\scriptstyle{g}}&&\rdTo^{\scriptstyle{h}}&
\\
X&&\rDashto^{\scriptstyle{\alpha}}&&Y
\\
\end{diagram}
\end{equation}
satisfying the following conditions:
\begin{enumerate}
\item
morphisms $g\colon W\to X$ and $h\colon W\to Y$ are contractions,
and the contraction $g$ is birational,
\item
if a prime divisor $A$ on $W$ is contracted by $g$, then it also
is contracted by $h$.
\end{enumerate}

\begin{example}
\label{ex-Pn-bundle}
Let $Y=\PP^1$ and let $W=\PP_Y (\sE)$ be the projectivization of
the vector bundle $\sE=\underbrace{\Oh_{\PP^1}\oplus\cdots
\oplus\Oh_{\PP^1}}_{n}\oplus\Oh_{\PP^1} (d)$, $d\ge 1$. Let
$g\colon W\to X\subset\PP^{d+n}$ be the morphism defined by the
linear system $\lins{\Oh_{\PP (\sE)} (1)}$. Then $g$ is a
birational contraction and we have diagram \eqref{diag-bss}. Here
$X\subset\PP^{d+n}$ is a cone over a rational curve
$C_d\subset\PP^{d}$ with vertex $\PP^{n-1}$ and the map $\alpha$
is defined by the linear system of planes on $X$. The base locus
of this system is precisely the vertex of the cone.
\end{example}

\begin{pusto}
For a rational $1$-contraction $\alpha\colon X\dasharrow Y$, we
may define the \emph{pull-back} of any $\RR$-Cartier divisor $D$
as follows: $\alpha^*D\eqdef g_*h^*D$ (it is easy to show that
this definition does not depend on the choice of the hut
\eqref{diag-bss}). Note however that the map $\alpha^*$ is not
functorial: it is possible that $(\alpha\comp\beta)^*$ does not
coincide with $\beta^*\alpha^*$. Similarly we can define
$\alpha_*$. For $\alpha_*$ we always have $\alpha_*\beta_*=
(\alpha\beta)_*$ whenever $\beta$ is a rational $1$-contraction.

For a closed subset $V\subset X$, denote by $ [V]^r$ the subset
consisting of all components $V_i\subset V$ of dimension $\ge r$.

A divisor $E$ on $X$ is said to be \emph{exceptional} with respect
to a rational $1$-contraction $\alpha$ if $\dim\alpha (E_i) <\dim
E_i$ and $\alpha (E_i)\neq Y$ for all components $E_i\subset\Supp
(E)$. In this situation $E$ is said to be \emph{very exceptional}
if for any prime divisor $G$ on $Y$ the divisorial fiber
$\alpha^\bullet (G)$ over $G$ is not contained in $\Supp (E)$ (see
\PLF{Def. 3.2} or \cite[Prop. 1.10] {F3}). Here the
\emph{divisorial fiber} is defined as $\alpha^\bullet (G)\eqdef
[\alpha^{-1} (G)]^{\dim X-1}$.
\end{pusto}

\begin{remark}
Let $\alpha$ be a morphism. To check that an exceptional divisor
$E$ is very exceptional it is sufficient to verify the following
property:
\begin{enumerate}
\item[]
for any component $E_i\subset\Supp (E)$ such that $\dim\alpha
(E_i)=\dim Y-1$ we have $\Supp (\alpha^*\alpha
(E_i))\not\subset\Supp (E)$.
\end{enumerate}

Assume that $\alpha$ is a morphism and let $d=\dim X -\dim Y$.
Define the following closed subset in $Y$:
\[
\mathfrak{E} (\alpha, E)\eqdef\left\{y\in Y\mid\left [\alpha\1
(y)\right]^{d}\subset\Supp (E)\right\},
\]
where $\alpha\1 (y) =f (h\1 (y))$. It is clear that $\mathfrak{E}
(\alpha, E)\subset\alpha (\Supp (E))$. Then $E$ is very
exceptional if and only if $\mt{codim}\mathfrak{E} (\alpha, E)\ge
2$.
\end{remark}

The following fact will be frequently used without references:

\begin{lemma}[{\cite[\S 1.1] {Sh}}, {\cite[2.15] {Sh1}}]
\label{11}
Let $f\colon X\to Z$ be a birational contraction and $A$ be an
$\RR$-Cartier divisor on $X$ such that
\begin{enumerate}
\item
$f$ contracts all components of $A$ with negative coefficients;
\item
for a sufficiently general curve $C$ in a fiber $A_i/f (A_i)$ in
each $f$-exceptional divisor $A_i$ having negative coefficient in
$A$, we have $A\cdot C\le 0$.
\end{enumerate}
Then the divisor $A$ is effective.
\end{lemma}

We need also the following result:

\begin{lemma}
\label{11a}
Let $f\colon X\to Z$ be a contraction \textup(not necessary
birational\textup) and let $A$ be a divisor on $X$. Write
$A=A^+-A^-$, where $A^+$ and $A^-$ are effective divisors without
common component. Assume that $-A$ is $f$-nef and $A^-$ is very
exceptional on $Z$. Then $A^-=0$, i.e., $A$ is effective.
\end{lemma}

\begin{proof}
Assume that $\dim X >\dim Z$. Let $A^-\neq 0$. If $\dim f (\Supp
(A^-)) > 0$, then we can replace $Z$ with its general hyperplane
section $Z'\subset Z$, $X$ with $f\1 (Z)'$, and $A$ with $A|_{f\1
(Z)'}$. The very exceptionality of $A^-$ is preserved. Indeed, it
is sufficient to choose a hyperplane section $Z'\subset Z$ so that
it does not contain components of the set $\mathfrak{E} (\alpha,
E)$ of codimension $2$. Continuing the process we get the
situation when $\dim f (\Supp (A^-)) =0$. We may also assume that
$Z$ is a sufficiently small affine neighborhood of some fixed
point $o\in Z$ (and $f (\Supp (A^-)) =o$). Further, all the
conditions of lemma are preserved if we replace $X$ with its
general hyperplane section $X'$. If $\dim Z > 1$, then we can
reduce our situation to the case $\dim X=\dim Z$. Then the
statement of the lemma follows by Lemma \xref{11} and from the
existence of the Stein factorization. Finally, consider the case
$\dim Z=1$ (here we may assume that $\dim X=2$). In this instance,
$A\cdot A^-=A^+\cdot A^-- (A^-)^2\le - (A^-)^2$. By the Zariski
lemma the last number is positive, a contradiction.
\end{proof}

\subsection*{Movable and fixed parts of a divisor}
For a divisor $D$, we put
\begin{equation}
\label{eq-Mov-Fix}\Mov (D)=
\begin{cases}
-\inf\limits_{s\in K (X)^*}\{(s)\mid D+ (s)\ge 0\} &\text{if
$|D|\neq\emptyset$},
\\
-\infty &\text{otherwise}.
\end{cases}
\end{equation}
If $H^0 (\Oh (D))\neq 0$ (that is equivalent to
$\lins{D}\neq\emptyset$), we put
\[
\Fix (D)=\inf\{L\mid L\sim D,\quad L\ge 0\}=D -\Mov (D).
\]
Divisors $\Mov (D)$ and $\Fix (D)$ are called \emph{mobile} and
\emph{fixed} parts of a divisor $D$ respectively. If $D$ is
effective, then so is $\Mov (D)$:
\[
\Mov (D)\ge - (\Const) =0.
\]
Obviously, $\Mov (D)=\Mov (\down {D})$.

\begin{lemma}[cf. {\PLF{Prop. 4.15}}]
\label{L-mobile-1}
Let $D$ be a divisor such that $H^0 (\Oh (D))\neq 0$. Then $M=\Mov
(D)$ satisfies the following properties:
\begin{enumerate}
\item
$M\le D$;
\item
the linear system $\lins{M}$ has no fixed components \textup(i.e.,
the divisor $M$ is b-free\textup);
\item
if a divisor $L\le D$ is b-free, then $L\le M$.
\end{enumerate}
Conversely, if an \textup(integral\textup) divisor $M$ satisfies
conditions {\rm (i)-(iii)}, then $M=\Mov (D)$.
\end{lemma}
\begin{proof}
It is clear that
\[
D-M=\inf_{s\in K (X)^*}\{D+ (s)\mid D+ (s)\ge 0\}\ge 0.
\]
Since
\[
0=\inf_{s\in K (X)^*}\{M+ (s)\mid D+ (s)\ge 0\}
\]
and $D\ge M$, we have that the divisor $D+ (s)$ is effective if
and only if so is $M+ (s)$. Hence,
\[
0=\inf_{s\in K (X)^*}\{M+ (s)\mid M+ (s)\ge 0\}.
\]
This means that the linear system has no fixed components.
Finally, let a divisor $L\le D$ b-free. Then
\[
L= -\inf_{s\in K (X)^*}\{(s)\mid L+ (s)\ge 0\}\le -\inf_{s\in K
(X)^*}\{(s)\mid D+ (s)\ge 0\}=M.
\]
The last statement follows from the uniqueness of a divisor $M$
satisfying conditions (i)-(iii).
\end{proof}

\begin{lemma}
\label{L-mobile-2}
Let $D$ be a divisor such that $H^0 (\Oh (D))\neq 0$. Then $M=\Mov
(D)$ satisfies the following properties:
\begin{enumerate}
\item
$M\le D$;
\item
$H^0 (\Oh (M)) =H^0 (\Oh (D))$;
\item
if for a divisor $L\le D$ the equality $H^0 (\Oh (L)) =H^0 (\Oh
(D))$ holds, then $L\ge M$.
\end{enumerate}
Conversely, if an \textup(integral\textup) divisor $M$ satisfies
conditions {\rm (i)-(iii)}, then $M=\Mov (D)$.
\end{lemma}

\begin{proof}
(ii) Since $M\le D$, $H^0 (\Oh (M))\subset H^0 (\Oh (D))$. Let
$s\in H^0 (\Oh (D))$. Then $D+ (s)\ge 0$ and according to
\eqref{eq-Mov-Fix} we have $M+ (s)\ge 0$, i.e., $s\in H^0 (\Oh
(M))$.

(iii) Let $H^0 (\Oh (L)) =H^0 (\Oh (D))$. Then $L+ (s)\ge 0$ if
and only if $D+ (s)\ge 0$. Hence,
\[
M= -\inf_{s\in K (X)^*}\{(s)\mid D+ (s)\ge 0\}= -\inf_{s\in K (X)
^*}\{(s)\mid L+ (s)\ge 0\}.
\]
Therefore,
\[
L-M=\inf_{s\in K (X)^*}\{L+ (s)\mid L+ (s)\ge 0\}\ge 0.
\]
\end{proof}

\section{Classical Zariski decomposition}
We say that a divisor $D$ is \emph{effective modulo $\QQ$-linear
\textup($\RR$-linear\textup) equivalence} if $D\qq D'$
(respectively, $D\qr D'$), where $D'$ is effective. The
effectiveness modulo $\QQ$-linear equivalence is equivalent to
that $H^0 (\Oh (\alpha D))\neq 0$ for some $\alpha\in\NN$.

Zariski decomposition and its various generalizations have the
form $D=P+N$, where $N$ is the effective part and $P$ is the
``maximal positive'' part. In general such a decomposition is
defined for effective divisors though many statements work in a
more general situation:

\begin{definition}[{\cite{F1}}]
\label{pseudoeff}
A divisor $D$ on projective variety $X$ is said to be
\emph{pseudo-effective} if there exists an ample divisor $H$ such
that $D+\ep H$ is effective for any $\ep > 0$.
\end{definition}

If $D$ is an $\RR$-Cartier divisor, then its pseudo-effectiveness
is equivalent to any of the following conditions (see \cite{Mo}):
\begin{enumerate}
\item
the class of $D$ is contained in the closure of the cone of
effective divisors (in numerical sense), i.e., there exists a
sequence of effective divisors $D^{(n)}$ such that $\lim\left
(D^{(n)}\cdot C\right) =D\cdot C$ for any curve $C$;
\item
the divisor $D+\ep H$ is effective (modulo $\qr$) for any ample
divisor $H$ and any $\ep > 0$;
\item
$D\cdot\ell\ge 0$ for any nef $1$-cycle $\ell$ on $X$.
\end{enumerate}

It is easy to see that the property of an $\RR$-Cartier divisor to
be pseudo-effective is closed under taking pull-backs $f^*$.

\begin{theorem}[\cite{Z}, \cite{F1}]
\label{main}
Let $X$ be a surface and let $D$ be a pseudo-effective divisor on
$X$. Then there exists an effective divisor $N=N (D)=\sum a_iN_i$
on $X$ such that
\begin{enumerate}
\item
the divisor $P\eqdef D-N$ is nef;
\item
either $N=0$ or the matrix $(N_i\cdot N_j)$ is negative definite;
\item
$(P\cdot N_i) =0$ for all $i$.
\end{enumerate}
Furthermore, if $D$ is a $\QQ$-divisor, then so is $N$. The
divisor $N (D)$ is uniquely defined by the class of numerical
equivalence of $D$.
\end{theorem}

A decomposition $D=N+P$ satisfying conditions (i) - (iii) of
Theorem~ \xref{main} is called a \emph{Zariski decomposition} of a
divisor $D$. The divisor $P$ is its \emph{positive} and $N$ is its
\emph{negative} (or \emph{exceptional} \PL) parts. Sometimes, by
abuse of language, we say that a Zariski decomposition is the
divisor $N=N (D)$.

Note that in the two-dimensional case the $\RR$- (or
$\QQ$)-Cartier condition is not necessary: the intersection theory
is defined for any normal surface (see, e.g.., \cite{Sa}). We give
a sketch of proof with running the ``$D$-MMP'' (see \cite{Sa}):

\begin{outline}
If $D$ is nef, then we put $N (D) =0$. Otherwise there exists an
irreducible curve $E$ such that $D\cdot E < 0$. Since $D$ is
pseudo-effective, we have $E^2 < 0$ (otherwise $E$ is a nef curve
and by definition \xref{pseudoeff}, $D\cdot E\ge 0$). By the
Grauert criterion, the curve $E$ is contractible (at least in the
category of normal analytic spaces): $f_1\colon X\to X_1$, where
the divisor $f_{1*}D$ is again pseudo-effective. Continuing the
process we obtain a model $X'$ on which the image $D'$ of $D$ is
nef (in other words, we run ``$D$-minimal model program''). Let
$f\colon X\to X'$ be the composition of all contractions. Put
$P\eqdef f^*D'$ and $N\eqdef D-P=D-f^*f_*D$. Since the divisor
$-D$ is nef over $X'$, we have $N\ge 0$ (see Lemma \xref{11}). The
uniqueness follows by Proposition \xref{prop-Zar-surf} below.
\end{outline}

Note however that in the category of projective surfaces
contraction $X\to X'$ does not necessarily exist. In other words,
the divisor $P$ is not always semiample. Our proof shows that we
may guarantee the semiampleness of $P$ for divisors of type
$D=K_X+B$ whenever the pair $(X, B)$ is log canonical.

In his paper \cite{Z} Zariski considered the case of effective and
integral divisor $D$. The generalization to the pseudo-effective
case belongs to Fujita \cite{F1}. Below, in Example~
\xref{ex-nonrat} we will see that in higher dimensional
generalizations one has to consider divisors with irrational
coefficients. We discuss the different higher dimensional
generalizations of Theorem \xref{main} and relations between them.

It is easy to see that Zariski decompositions agree via
pull-backs. That is why we have the following.

\begin{proposition}
For any pseudo-effective divisor $D$ on a surface $X$, there
exists a decomposition $\ov D=\sP+\sN$ of b-divisors in the group
$\bCDiv_{\RR} (K (X))$ of b-Cartier b-divisors \textup
(see\Isk\textup) that induce Zariski decompositions $(\ov D)
_Y=\sP_Y+\sN_Y$ on all normal projective models $Y/X$.
\end{proposition}
Simple examples show also that Zariski decompositions does not
agree via $f_*$. That is why $\sP_Y$ and $\sN_Y$ are not elements
of $\bDiv_{\RR} (K (X))$.

\begin{example}
\label{ex-cremona}
Consider the quadratic birational map $\PP^2\dashrightarrow\PP^2$,
$(x: y: z)\dashrightarrow (yz: xz: xy)$ and let
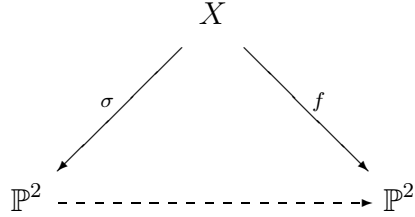
\begin{equation*}
\begin{diagram}
&&X&&
\\
&\ldTo^{\scriptstyle{\sigma}}&&\rdTo^{\scriptstyle{f}}&
\\
\PP^2&&\rDashto&&\PP^2
\\
\end{diagram}
\end{equation*}
be a resolution of the indeterminacy. Here $\sigma\colon X\to
\PP^2$ and $f\colon X\to\PP^2$ are blowups of triples of different
points on $\PP^2$. Let $D\eqdef\sigma^*H+E_1+E_2+E_3$, where $H$
is an ample divisor on $\PP^2$ and let $E_1, E_2, E_3$ be
exceptional divisors of $\sigma$. Then $N (D) =E_1+E_2+E_3$ and $P
(D)=\sigma^*H$. But $N (f_*D) =0\neq f_* (E_1+E_2+E_3)$ and $P
(f_*D) =f_*D\neq f_*\sigma^*H$.
\end{example}

\begin{proposition}
\label{prop-Zar-surf}
Let $X$ be a surface and let $D$ be a pseudo-effective divisor on
$X$. Let $D=N (D) +P (D)$ be a Zariski decomposition. Then for any
nef divisor $L$ such that $L\le D$ we have $L\le P (D)$.
\end{proposition}

\begin{proof}
Write $N=N (D)=\sum a_iN_i$. Let $L\le D$, where $L$ is nef. For
all $i$ we have
\[
N_i\cdot (P-L)=- N_i\cdot L\le 0.
\]
Write $P-L= (D-L)-N=F^{\sharp}-N^{\sharp}$, where $F^{\sharp}$ and
$N^{\sharp}$ are effective divisors without common components.
Since $N^{\sharp}\le N$,
\[
0\ge N^{\sharp}\cdot (P-L) =N^{\sharp}\cdot F^{\sharp}-N^{\sharp
2}\ge 0.
\]
This gives us $N^{\sharp 2}=0$, $N^{\sharp}=0$ and $P-L\ge N$.
\end{proof}
Proposition \xref{prop-Zar-surf} follows also by Lemma \xref{11}
applied to the contraction $f$ from the proof of Theorem
\xref{main}.

\begin{pusto}\textbf{Generalizations.}
\label{gen-prin}
It is clear that Zariski decomposition cannot be generalized in
higher dimensions without significant modifications. We formulate
general scheme for eventual generalizations. Let $D$ be an
$\RR$-Cartier divisor on a variety $X$ and let $N (D)$ be the
negative part in a (generalized) Zariski decomposition. It is
reasonable to claim the following:
\begin{enumerate}
\item
the divisor $N (D)$ is effective;
\item
the divisor $D-N (D)$ is ``positive'' in some sense;
\item
$N (D)$ is ``minimal''.
\end{enumerate}
Also it is reasonable to claim that the decompositions agree via
pull-backs: $N (D)$ is $\RR$-Cartier and if $f\colon Y\to X$ is a
birational contraction, then
\begin{equation}
\label{surd}
N (f^*D) =f^*N (D).
\end{equation}
\end{pusto}

\begin{remark}
\label{zam-b-div}
The last property allow us to introduce Zariski decompositions of
b-Cartier b-divisors: for $\sD\in\bCDiv_{\RR} (K (X))$ there
exists $\sN (\sD)\in\bCDiv_{\RR} (K (X))$ and a projective model
$Y$ of field $K (X)$ such that $\sN (\sD) _{Y'}$ is a
(generalized) Zariski decomposition for any model of $Y'$ of $K
(X)$ dominating $Y$.
\end{remark}

\section{On finite generation of divisorial algebras}
One of fundamental problems of algebraic geometry is the question
about finite generation of algebras
\[
\sR_X D=\bigoplus_{n\ge 0} H^0 (X,\Oh_X (nD)).
\]
Here $\sR_X D$ is considered as a subalgebra of the algebra $K (X)
[t]$, where each space $H^0 (\Oh_X (nD))$ is enclosed in the
component $K (X) t^n$. Below we present several well known facts.

\begin{proposition}[see, for example, {\cite[ch. III, \S 1, 3${}^o$]
{Bu}}, {\PLF{Th. 4.6}}]
\label{prop-truncation}
Let $\sR=\bigoplus_{n\ge 0}\sR_n$ be a graded algebra over a field
$k$ such that $\sR_0=k$. Then
\begin{enumerate}
\item
the algebra $\sR$ is finitely generated if and only so is the
\emph{truncated} algebra $\sR^{[n_0]}\eqdef\bigoplus_{n\ge
0}\sR_{nn_0}$, $n_0\in\NN$;
\item
if the algebra $\sR$ is finitely generated, then there exists
$n_0\in\NN$ such that the truncated algebra $\sR^{[n_0]}$ is
generated by elements degrees $1$.
\end{enumerate}
\end{proposition}

\begin{proposition}
Let $D$ be an integral divisor such that $\Bs\lins{D}=\emptyset$.
Then the algebra $\sR_X D$ is finitely generated.
\end{proposition}
\begin{proof}
First, we prove this in the case, when $D$ is ample. Let
$X\subset\PP^N$ be an embedding corresponding to a suitable
multiplicity $n_0D$ of $D$ and let
$\mathcal{J}_X\subset\Oh_{\PP^N}$ be the ideal sheaf. From the
exact sequence
\[
0\longrightarrow\mathcal{J}_X (n)\longrightarrow\Oh_{\PP^N}
(n)\longrightarrow\Oh_X (nn_0D)\longrightarrow0
\]
and Serre's vanishing theorem we obtain that the restrictions $H^0
(\Oh_{\PP^N} (n))\longrightarrow H^0 (\Oh_X (n_0D))$ are
surjective for $n\ge n_1$. Hence there is a surjective map
$\sR_{\PP^N}^{[n_1]}\Oh (1)\to\sR_{X}^{[n_0n_1]} D$. According to
Proposition \xref{prop-truncation} this is proves finite
generation of $\sR_X D$.

Now consider the general case. Let $X\to\bar X\subset\PP^N$ be a
morphism defined by the linear system $|n_0D|$, where $n_0\gg 0$
and $\bar X=\var (X)$. Consider its Stein factorization
$X\stackrel{\var}{\longrightarrow}
X'\stackrel{\psi}{\longrightarrow}\bar X$. Then the divisor
$H=\psi^*\Oh_{\bar X} (1)$ is ample and $\Mov (nD)=\var^*H$. Since
$\var_*\Oh_X=\Oh_{X'}$, $\sR^{[n]}D\simeq\sR_{X'} H$. According to
the above the last algebra is finitely generated.
\end{proof}

\subsection*{Zariski decomposition and finite generation of
divisorial algebras}

\begin{pusto}
For a Zariski decomposition $D=P+N$, there is an isomorphism of
graded algebras
\begin{equation}
\label{eq-R-sim-R}
\sR_X D\simeq\sR_X P.
\end{equation}
Indeed, for each $n\in\NN$ we have
\[
H^0 (\Oh_X (nP))\subset H^0 (\Oh_X (nD)) =H^0 (\Oh_X (\Mov (nD))).
\]
On the other hand, the divisor $\Mov (nD)$ is nef and by
Proposition \xref{prop-Zar-surf} we have $\Mov (nD)\le nP$. This
gives us the inverse inclusion
\[
H^0 (\Oh_X (\Mov (nD)))\subset H^0 (\Oh_X (nP)).
\]
\end{pusto}

Thus the question about finite generation of a divisorial algebra
$\sR_X D$ can be reduced to the question about finite generation
of the divisorial algebra $\sR_X P$, where the divisor $P$ is nef.
It is well known that the algebra $\sR_X D$ is not always finitely
generated:

\begin{proposition}[\cite{Z}, cf. {\PLF{Th. 4.28}}]
\label{ex-Zar-prop}
Suppose that an integral effective Cartier divisor $D$ satisfies
the following two conditions:
\begin{enumerate}
\item
$\Fix\lins{nD}\neq\emptyset$ for all $n\in\NN$;
\item
multiplicities of components $\Fix\lins{nD}$ are bounded as
$n\to\infty$.
\end{enumerate}
Then the algebra $\sR_X D$ is not finitely generated.
\end{proposition}

The proposition above enable us to construct a great number of
divisors with non-finitely generated algebras $\sR_X D$ (see
Example \xref{ex-ell-Zar} below).

\begin{proof}
Suppose that the algebra $\sR_X D$ is generated by the finite
number of elements $u_1,\dots, u_{r}$. Let $d_i=\deg u_i$ and
$d\eqdef\max\{d_1,\dots, d_r\}$. Then the vector space $H^0 (\Oh_X
(nD))$ is generated by the monomials of the form
\[
u_{1}^{\nu_{1}}u_{2}^{\nu_{2}}\cdots u_{r}^{\nu_{r}},
\]
where
\[
\nu_i\in\ZZ_{\ge 0},\qquad\sum_{i=1}^{r} d_i\nu_{i}=n.
\]
Therefore,
\[
\Fix\lins{nD}\ge\Min\left\{\left.\sum_{i=1}^{r}\nu_i\Fix\lins{d_iD}\
\right|\ \sum_{i=1}^{r} d_i\nu_{i}=n\right\}.
\]
On the other hand,
\[
\frac{d!}{k}\Fix\lins{kD}\ge\Fix\lins{d! D}\neq 0
\]
for all $k=1,\dots, d$. Thus divisors $\Fix\lins{d_1D}$,\dots,
$\Fix\lins{d_rD}$ have at least one common component. This gives
us a contradiction.
\end{proof}

The condition (ii) in Proposition \xref{ex-Zar-prop} is
automatically satisfied if the divisor $D$ is nef and big:

\begin{corollary}
\label{ex-Zar-prop-11}
Let $D$ be a nef and big integral Cartier divisor on a surface $X$
such that $\Fix\lins{nD}\neq\emptyset$ for all $n\in\NN$. Then the
algebra $\sR_X D$ is not finitely generated.
\end{corollary}

\begin{proof}
Follows by Proposition \xref{prop-C-W} below.
\end{proof}

We say that the base locus $\Bs\lins{nD}$ is \emph{bounded as
$n\to\infty$} if for any birational contraction $f\colon Y\to X$
multiplicities of components of $\Fix\lins{nf^*D}$ are bounded as
$n\to\infty$. In other words multiplicities of components of the
b-divisor $\Fix (\ov{nD})$ are bounded (see
\xref{def-pseudo-Mov-Fix}).

\begin{proposition}[{\cite[Th. 2.2] {W}}]
\label{prop-C-W}
Let $D$ be an \textup(integral\textup) Cartier divisor on a
variety $X$. If the base locus $\Bs\lins{nD}$ is bounded as
$n\to\infty$, then the divisor $D$ is nef. Conversely, if $D$ is
nef and big, then $\Bs\lins{nD}$ is bounded.
\end{proposition}
\begin{outline}
We give an outline of the proof only for $\dim X=2$. If $D\cdot C
< 0$ for some curve $C$, then $C\subset\Fix\lins{nD}$ for any
$n\in\NN$. Furthermore, $(nD-mC)\cdot C < 0$ whenever $m <
n\frac{D\cdot C}{C^2}$, i.e., $mC\le\Fix\lins{nD}$ for such $m$.
This is means that the multiplicity $C$ in $\Fix\lins{nD}$ goes to
infinity, a contradiction.

Conversely, assume that $D$ is nef and big. Let $E_i$ be fixed
components of $\lins{nD}$. Obviously, we may assume that the
surface $X$ and all the $E_i$ are nonsingular. Choose a very ample
divisor $H$ such that divisors $H-K_X-E_i$ are ample for all
$E_i$. We prove that $H^0 (E_i,\OOO_{E_i} (mD+H))\neq 0$. By the
Riemann-Roch Theorem this is satisfied if
\[
(mD+H)\cdot E_i > p_a (E_i)-1=\frac12 (K_X+E_i)\cdot E_i.
\]
The last is equivalent to
\[
0 < (2mD+2H-K_X-E_i)\cdot E_i=2mD\cdot E_i+H\cdot E_i+
(H-K_X-E_i)\cdot E_i,
\]
that, obviously, is satisfied. Further, $H^1 (\OOO_X (mD+H-E_i))
=0$. From exact sequence
\begin{multline*}
0\longrightarrow H^0 (\OOO_X (mD+H-E_i))\longrightarrow H^0
(\OOO_X (mD+H))
\\
\longrightarrow H^0 (\OOO_{E_i} (mD+H))\longrightarrow 0
\end{multline*}
we obtain that $E_i$ is not a fixed component of $\lins{mD+H}$.
Therefore the linear system $\lins{mD+H}$ has no fixed component
for $m\in\ZZ_{\ge 0}$. Finally, since $D$ is big, for some
$a\in\NN$ we have $aD\sim H+F$, where $F$ is an effective divisor.
Therefore, $\Fix\lins{(a+m) D}=\Fix\lins{F+H+mD}$ is bounded by
the divisor $F$, a contradiction.
\end{outline}

From Corollary \xref{ex-Zar-prop-11} (see also Example
\xref{ex-D=0,D=1} below) taking into account isomorphism
\eqref{eq-R-sim-R} we obtain the following criterion.

\begin{theorem}[\cite{Z}, cf. {\PLF{Th. 4.28}} and
Theorem {\xref{th-Lim-Crit}}] Let $X$ be a surface and let $D$ be
an effective modulo $\QQ$-linear equivalence $\QQ$-Cartier divisor
on $X$. Then the algebra $\sR_X D$ is not finitely generated if
and only if $\kappa (X, D) =2$ and the divisor $P (D)$ is not
semiample.
\end{theorem}

\begin{example}[Zariski]
\label{ex-ell-Zar}
Consider a nonsingular cubic curve $C\subset\PP^2$. Pick 12 points
$P_1,\dots, P_{12}\in C$ so that the divisor $\Oh_C (4) -\sum P_i$
is not a torsion in $\Pic (C)$. Let $\sigma\colon X\to\PP^2$ be
the blowup of $P_1,\dots, P_{12}$ and let $E_1,\dots, E_{12}$ be
the corresponding exceptional divisors. Put
$D\eqdef\sigma^*\Oh_{\PP^2} (4) -\sum E_i$. It is easy to show
that the divisor $D$ is effective modulo linear equivalence, nef
and big. Therefore in a Zariski decomposition one has $N (D) =0$.
Furthermore, $D\cdot\widetilde C=0$ and the birational transform
$\widetilde C\eqdef\sigma\1 (C)$ of the curve $C$ is a fixed
component of the linear system $\lins{nD}$ for any $n\in\NN$
(otherwise $nD|_{\widetilde C}=0$). Therefore $D$ satisfies the
conditions of Proposition \xref{ex-Zar-prop} and the algebra
$\sR_X D$ is not finitely generated.
\end{example}

For a big divisor, there is also the following criterion of finite
generation of the algebra $\sR_X D$:
\begin{theorem}[{\cite[Th. 1.2] {W}}, cf. {\PLF{Th. 3.18}}]
\label{prop-C-W-finite}
Let $D$ be an \textup(integral\textup) big Cartier divisor. Then
the following conditions are equivalent:
\begin{enumerate}
\item
the algebra $\sR_X D$ is finitely generated;
\item
there exists $n\in\NN$ and a birational contraction $f\colon\tilde
X\to X$ such that the linear system $\Mov\lins{nf^*D}$ defines a
birational morphism contracting all the components of
$\Fix\lins{nf^*D}$.
\end{enumerate}
\end{theorem}

\section{s-, $\sigma$- and sectional decompositions}

\subsection*{s-decomposition}
Let $D$ be an effective modulo $\QQ$-linear equivalence divisor.
Put
\[
M_n\eqdef\Mov (nD).
\]
Thus $nD=M_n+F_n$, where $F_n\ge 0$. Put also
\begin{equation}
\label{eq-def-s} P_{\s}
(D)\eqdef\limsup_{n\to\infty} (M_n/n).
\end{equation}
Since $M_{n_1n_2}\ge n_2M_{n_1}$ for all $n_1, n_2\in\NN$, we have
\[
nP_{\s} (D)\ge M_n\quad\text{for all $n\in\NN$.}
\]
We obtain a decomposition
\[
D=P_{\s} (D) +N_{\s} (D),\qquad N_{\s} (D)\ge 0,
\]
which we call an \emph{s-decomposition}. The divisor $P_{\s} (D)$
is called its \emph{positive} and $N_{\s} (D)$ its \emph{negative}
part. In case, when $D$ is not effective modulo $\QQ$-linear
equivalence, we put $P_{\s} (D)=-\infty$. Obviously,
\[
N_{\s} (D)=\inf_{
\begin{array}{c}
\scriptscriptstyle{n\in\NN}
\\
\scriptscriptstyle{s\in K (X)^*}
\end{array}}\{D+ (s)/n\mid D+ (s)/n\ge 0\}.
\]
In other words
\begin{equation}
\label{eq-Ns-}
N_{\s} (D)=\inf\{L\mid L\qq D,\ L\ge 0\}.
\end{equation}

Decompositions of such type were used by many authors: \cite[\S
2]{Kaw-crep}, \cite[\S 2]{Nak}, \PLF{\S\S 3-4}. Note that $P_{\s}
(D)$ and $N_{\s} (D)$ are not necessarily $\QQ$-divisors, even if
the divisor $D$ is integral (see Example \xref{ex-nonrat}).

\begin{remark}
\label{rem-s-Zar}
Let $X$ be a surface and let $D$ be an effective modulo $\qq$
divisor on $X$. Then the divisor $P_{\s}$ is nef and by
Proposition \xref{prop-Zar-surf} we have $P\ge P_{\s}$, where $P=P
(D)$ is the positive part of the classical Zariski decomposition.
In general case, this inequality is not an equality (see Example
\xref{ex-9-points}). If $\kappa (X, D)\ge 1$ and surface
$\QQ$-factorial, then equality $P= P_{\s}$ holds, i.e., the
s-decomposition coincides with the classical Zariski
decomposition.

Indeed, if the divisor $D$ is big, then so is the divisor $P$ and
by the Kodaira lemma $P=A+F$, where $A$ is ample and $F$ is an
effective divisor. Since for any $0 <\ep < 1$, the divisor $(1
-\ep) P+\ep A$ is ample and $(1 -\ep) P+\ep A\le P$, we have $(1
-\ep) P+\ep A\le P_{\s}$. Therefore, $P\le P_{\s}$ (see also
Proposition \xref{sect-dec-ed} below).

If $\kappa (X, D)=1$, then $\kappa (X, P)=1$, $P^2=0$ and for some
$n\in\NN$ the linear system $\Mov (nP)$ defines a contraction
$f\colon X\to Z$ onto a curve. It is easy to see that
\begin{multline*}
P\cdot\Mov (nP) =P\cdot\Fix (nP)=\Mov (nP)^2=
\\
\Mov (nP)\cdot\Fix (nP)=\Fix (nP)^2=0.
\end{multline*}
Hence the divisor $\Fix (nP)$ is contained in fibers and is a
pull-back of some divisor on $Z$. Therefore we can write $nP\sim
f^*L$. Then $N_{\s} (P)=\frac1nN_{\s} (f^*L)$ (see \xref{svoistva}
below). Since $\Mov (mf^*L)\ge f^*\Mov (mL)$, $P_{\s} (f^*L)\ge
f^*P_{\s} (L) =f^*L$. Thus $P_{\s} (f^*L) =f^*L$, $N_{\s} (P)
=N_{\s} (f^*L) =0$, $P_{\s} (P) =P$ and $P_{\s} (D)\ge P_{\s} (P)
=P$.
\end{remark}

\subsection*{$\sigma$-decomposition}
In the work \cite{Nak} Nakayama defined a similar type of
decompositions for any pseudo-effective divisor $D$, so-called,
\emph{$\sigma$-decomposition}:
\[
D=P_{\sigma} (D) +N_{\sigma} (D).
\]
If $D$ is a big divisor, then
\[
N_{\sigma} (D)\eqdef N_{\s} (D)=\inf\{L\mid L\qq D,\quad L\ge 0\}.
\]
If $D$ is a pseudo-effective, but not big divisor, then we put
\[
N_{\sigma} (D)\eqdef\lim_{\ep\to 0} N_{\sigma} (D+\ep A),
\]
where $A$ is an arbitrary ample divisor. (It is easy to show that
this definition does not depend on the choice of $A$).

From now on $\ov\Mv (X)$ denotes the closed convex cone in $\N^1
(X)$ generated by the classes of mobile Cartier divisors and
$\Mv^o (X)$ denotes the interior of $\ov\Mv (X)$. By the Kodaira
lemma ewe have that if the class of any $\QQ$-divisor $D$ is
contained in $\Mv^o (X)$, then some multiplicity $nD$, $n\in\NN$
is an integral mobile divisor \cite[\S 2]{Kaw-crep}.

\begin{proposition}[{\cite[2.1. 10] {Nak}}]
\label{prop-second-s-decomp}
$\sigma$-decomposition of a pseudo-effective divisor $D$ on a
nonsingular \textup(enough: on $\QQ$-factorial\textup) variety
satisfies the following properties:
\begin{enumerate}
\item
$ [P_\sigma (D)]\in\ov{\Mv} (X)$;
\item
for any divisor $L$ such that $L\le D$ and $[L]\in\ov{\Mv} (X)$ we
have $L\le P_\sigma (D)$.
\end{enumerate}
In particular, on a nonsingular surface the $\sigma$-decomposition
coincides with the classical Zariski decomposition.
\end{proposition}

Thus $P_{\sigma} (D)\ge P_{\s} (D)$ for any effective modulo $\qq$
divisor $D$. If the divisor $D$ is big, then (by definition)
$P_{\sigma} (D) =P_{\s} (D)$. However, this is not true for
arbitrary effective divisors even in the two-dimensional case (see
Example \xref{ex-9-points} below).

\begin{example}
\label{ex-9-points}
Let $C\subset\PP^2$ be a nonsingular cubic curve and let
$P_1,\dots, P_9\in C$ be distinct points such that $\Oh_C (3)
-\sum P_i$ is not a torsion in $\Pic (C)$. Let $\sigma\colon
X\to\PP^2$ be the blowup of points $P_1,\dots, P_9$ and let $D$ be
the birational transform $C$. Then $\dim\lins{nD}=0$ for all
$n\in\NN$. Therefore, $P_{\s} (D) =0$. On the other hand, $D$ is
nef. Hence, $P_{\sigma} (D) =P (D) =D$.
\end{example}

\subsection*{Properties of s-decompositions}
The following properties are immediate consequences of the
definition.
\begin{pusto}
\label{svoistva}
\begin{enumerate}
\item
The negative part of an s-decomposition depends only on the class
of $\QQ$-linear equivalence of the divisor $D$. More precisely, if
$D$ and $D'$ are $\QQ$-linearly equivalent and effective modulo
$\qq$ divisors, then $N_{\s} (D) =N_{\s} (D)'$.
\item
If $D'$ and $D''$ are effective modulo $\qq$ divisors, then
\[
P_{\s} (D)' +P_{\s} (D'')\le P_{\s} (D'+D''),\qquad N_{\s} (D)'
+N_{\s} (D'')\ge N_{\s} (D'+D'').
\]
If additionally $D\ge D'$, then $P_{\s} (D)\ge P_{\s} (D)'$.
\item
Let $D$ be an effective modulo $\qq$ divisor. Then for any
$\alpha\in\QQ_{> 0}$ we have $N_{\s} (\alpha D)=\alpha N_{\s}
(D)$.
\item
Computing $P_{\s}$ we always can replace the limit
\eqref{eq-def-s} on the ``truncated'' limit:
\[
P_{\s} (D)=\limsup_{n\to\infty}\frac{M_{nn_0}}{nn_0}.
\]

\item
For an ample divisor $D$, we have $P_s (D) =D$.
\end{enumerate}
\end{pusto}

\begin{proposition}
\label{prop-s-dec-semiample}
Let $D$ be an effective modulo $\QQ$-linear equivalence divisor.
If a $\QQ$-divisor $L$ is b-semiample and $L\le D$, then $L\le
P_{\s} (D)$.
\end{proposition}

\begin{proof}
Write $L=\sum\alpha_i H_i$, where $H_i$ are integral b-free
divisors and $\alpha_i\in\RR_{\ge 0}$. Since the coefficients of
the divisor $L$ are rational, we can choose numbers $\alpha_i$
also to be rational. Therefore $nL$ is an (integral) b-free
divisor for some $n\in\NN$. By Lemma \xref{L-mobile-1} we have
$L\le M_n/n\le P_{\s}$.
\end{proof}

The following easy statement shows how s-decompositions can be
used in the study of divisorial algebras.
\begin{proposition}
\label{prop-first-s-decomp}
An s-decomposition of an effective modulo $\QQ$-linear equivalence
divisor satisfies the following properties:
\begin{enumerate}
\item
$P_{\s} (D)\le D$;
\item
$\sR_XD=\sR_XP_{\s} (D)$;
\item
for any divisor $L$ such that $L\le D$ and $\sR_XD=\sR_XL$ we have
$L\ge P_{\s} (D)$ \textup(i.e., $P_{\s}$ is the \emph{smallest}
divisor satisfying properties {\rm (i)} and {\rm (ii)}\textup).
\end{enumerate}
\end{proposition}

Thus this proposition and \PLF{Remark 3.30} explain and justify
introduction of the concept ``s-decomposition''.

\begin{proof}
Since $P_{\s}\ge M_n/n$ we have by Lemma \xref{L-mobile-2} that
for any $n\in\NN$ we have
\[
H^0 (\Oh_X (M_n)) =H^0 (\Oh_X (nD))\supset H^0 (\Oh_X
(nP_{\s}))\supset H^0 (\Oh_X (M_n)).
\]
This proves (ii).

We prove (iii). Since $H^0 (\Oh (nD)) =H^0 (\Oh (nL))$, $nL\ge
M_n$ (again by Lemma \xref{L-mobile-2}). Hence, $L\ge P_{\s} (D)$.

The converse follows from the fact that there exists at most one
divisor satisfying properties (i) - (iii).
\end{proof}

\begin{remark}
Proposition \xref{prop-first-s-decomp} remains to be true if we
replace the condition (iii) with the following:
\begin{enumerate}
\item[(iii${}'$)]
for any divisor $L$ such that $L\le P_{\s} (D)$ and
$\sR_XD=\sR_XL$ we have $L=P_{\s} (D)$.
\end{enumerate}
(i.e., $P_{\s}$ is the \emph{minimal} divisor satisfying
properties {\rm (i)} and {\rm (ii)}). Indeed, by Proposition
\xref{prop-first-s-decomp} the s-decomposition satisfies
properties (i), (ii), (iii)${}'$. Conversely, let $P_{\s}'$
satisfies properties (i), (ii), (iii)${}'$. Then from (iii)
applied to $P_{\s}$ we have $P_{\s}'\ge P_{\s}$ and from
(iii)${}'$ we obtain $P_{\s}'=P_{\s}$.
\end{remark}

\begin{theorem}[cf. \PLF{Th. 4.28}]
\label{th-Lim-Crit}
Let $D$ be an effective modulo $\QQ$-linear equivalence divisor.
If the divisorial algebra $\sR_XD$ is finitely generated, then
$P_{\s} (D)=\Mov (n_0D)/n_0$ for some $n_0\in\NN$ \textup(in other
words the limit \eqref{eq-def-s} \emph{stabilizes}\textup).
\end{theorem}

Thus in the case when the algebra $\sR_XD$ is finitely generated,
the positive part of s-decomposition is a b-semiample
$\QQ$-divisor. Under additional conditions for the positive part
we obtain decompositions discussed in \S\S
\xref{sect-Shok}-\xref{sect-Fudj}. Theorem \xref{th-Lim-Crit} is
not a criterion: in this form the converse is not true. In fact,
the condition $P_{\s} (D)=\Mov (n_0D)/n_0$ is divisorial and is
preserved under small birational contractions, while finite
generation is essentially more subtle condition. In order to
obtain a criterion of finite generation, we must consider the
condition for the stabilization of limits on all blowups of the
initial variety, i.e., to pass to a b-divisor (see \PLF{Th. 4.28}
and Theorem \xref{th-Lim-Crit-b-div}). We notice that the
condition of stabilization does not mean that $P_{\s} (D)=\Mov
(nD)/n$ for all $n\gg 0$. However, this condition implies that
$P_{\s} (D)=\Mov (nn_0D)/(nn_0)$ for all $n\in\NN$.

\begin{proof}
According to Proposition \xref{prop-truncation} there exists
$n_0\in\NN$ such that the algebra $\sR_X (n_0D)$ is generated by
the elements $u_1,\dots, u_r\in H^0 (X,\OOO_X (n_0D))^*=H^0
(X,\OOO_X (n_0D))\setminus\{0\}$. For any $n\in\NN$ and for any
$s\in H^0 (X,\OOO_X (nn_0D))$, we have
\[
s=\sum_{\nu_1+\cdots+\nu_r=n} a_{\nu_1,\dots,\nu_r}
u_1^{\nu_1}\cdots u_r^{\nu_r},\quad\nu_1,\dots,\nu_r\in\ZZ_{>
0},\quad a_{\nu_1,\dots,\nu_r}\in\CC.
\]
Hence,
\begin{multline*}
\ord_F (s)\ge\min_{\nu_1+\cdots+\nu_r=n} (\nu_1\ord_F
(u_1)+\cdots+\nu_r\ord_F (u_r))
\\
\ge n\inf_{u\in H^0 (\OOO_X (n_0D))^*}\ord_F (u).
\end{multline*}
By definition we obtain $\Mov (nn_0D)\le n\Mov (n_0D)$. Since the
inverse inequality always holds, we have $\Mov (nn_0D)=n\Mov
(n_0D)$. This proves the theorem.
\end{proof}

\begin{examples}
\label{ex-D=0,D=1}
(i) If $\kappa (X, D) =0$, then $M_n=0$ or $-\infty$ for all $n$.
Moreover there exists $n_0\in\NN$ such that $H^0 (\Oh (n_0mD)) =1$
for all $m\in\NN$. Therefore, $M_{n_0m}=0$. This gives us $P_{\s}
(D) =0$.

(ii) If the class of the divisor $D$ is contained in the open cone
$\Mv^o (X)$, then $D$ we may be approximated from below by
b-semiample $\QQ$-divisors. Therefore, $P_{\s} (D) =D$.

(iii) Let $D$ be an effective divisor with $\kappa (X, D) =1$.
Then there exists $n_0\in\NN$ such that $\dim\lins{nD} > 0$ for
all $n\ge n_0$. Thus for every $n\ge n_0$ the linear system
$\lins{M_n}$ defines a map $X\dashrightarrow Y_n\subset\PP^{N_n}$.
It is clear that $K (Y_n)$ is a subfield in $K (X)$ generated by
$H^0 (\Oh_X (M_n))$. Since $M_n\le M_{n+1}$, $K (Y_n)\subset K
(Y_{n+1})$. All the fields $K (Y_n)$ have transcendence degree $1$
over $\CC$. Therefore the algebraic closure $\overline{K
(Y_{n_0})}$ of $K (Y_{n_0})$ in $K (Y)$ also has transcendence
degree $1$. Then there exists $n_1\ge n_0$ such that $K
(Y_n)=\overline{K (Y_{n_0})}$ for $n\ge n_1$. Let $Y$ be a
nonsingular curve such that $K (Y)=\overline{K (Y_{n_0})}$. All
rational maps $X\dashrightarrow Y_n$, $n\ge n_0$ factorise through
$Y$:
\[
X\stackrel{g}{\dashrightarrow} Y\stackrel{h_n}{\longrightarrow}
Y_n,
\]
where $h_n$ is a finite morphism. Therefore,
\[
M_n=g^*\bigl (\sup\{L\in\Div (Y)\mid nD-f^*L\ge0\}\bigr)
\]
for $n\gg 0$. Hence the divisor
\[
P_{\s} (D)=\lim_{n\to\infty} M_n/ n=g^*\bigl (\sup\{F\in\Div_{\QQ}
(Y)\mid D-f^*F\ge0\}\bigr)
\]
is b-semiample. Obviously, the last divisor is a $\QQ$-divisor
whenever so is $D$. In particular, we obtain that if $D$ is a
$\QQ$-divisor, then the algebra $\sR_XD\simeq\sR_XP_{\s}$ is
finitely generated. The constructed map is a particular case of
the \emph{Iitaka fibration} of $(X, D)$.

(iv) (see also Proposition \xref{prop-C-W}) Suppose that an
$\RR$-Cartier divisor $D$ is nef and big. Then by the Kodaira
lemma there exists an effective divisor $F$ such that the divisor
$D-F$ is ample. Therefore so is $D -\ep F$ for any $0 <\ep < 1$.
Thus,
\[
P_{\s} (D)\ge P_{\s} (D -\ep F) =D -\ep F.
\]
Passing to the limit we obtain $P_{\s} (D) =D$ (cf. with Example
\xref{ex-9-points}).
\end{examples}

\begin{remark}
If $N_{\s} (D) =0$, then the divisor $D$ is the limit of
b-semiample divisors $M_n/n$, i.e., its class is contained in the
closure of the cone $\Mv (X)$. In general case, this does not
imply the b-semiampleness of $D$ (see Example \xref{ex-ell-Zar}).
\end{remark}

\begin{proposition}[cf. {\cite{Ka}}]
\label{prop-s-decomp-unique}
Let $D$ be a big divisor on a $\QQ$-factorial variety. Then the
negative part of s-decomposition depends only on the class of
numerical equivalence of $D$.
\end{proposition}

\begin{proof}
Let $D$, $D'$ be big numerically equivalent divisors and let
$T\eqdef D'-D\equiv 0$. Since the divisor $P_{\s} (D)$ is big, the
Kodaira lemma gives us $P_{\s} (D) =H+F$, where $F\ge 0$ and $H$
is an ample divisor. Write
\[
D'=D+T=P_{\s} (D) +T+N_{\s} (D)=P_{\s} (D) -\ep F+T+\ep F+N_{\s}
(D).
\]
For any $0 <\ep < 1$, the class of $P_{\s} (D) -\ep F+T= (1 -\ep)
P_{\s} (D)+\ep H+T$ is contained in the cone $\Mv^o (X)$.
According to Example \xref{ex-D=0,D=1}, (ii) we have
\[
P_{\s} (D)'\ge P_{\s} (P_{\s} (D) -\ep F+T) =P_{\s} (D) -\ep F+T.
\]
Passing to the limit as $\ep\to 0$, we obtain $P_{\s} (D) +T\le
P_{\s} (D)'$, i.e., $N_{\s} (D)'\le N_{\s} (D)$. By symmetry we
have also the inverse inequality.
\end{proof}

\begin{proposition}
Let $f\colon X\to Y$ be a birational contraction and let $D$ be an
effective modulo $\qq$ divisor on $X$. Then

\begin{equation}
\label{eq-push-down}
f_* P_{\s} (D)\le P_{\s} (f_*D).
\end{equation}
\end{proposition}
\begin{proof}
For any $n\in\NN$, we have $f_*\Mov (nD)\le f_*nD$ and the divisor
$f_*\Mov (nD)$ is b-free. Hence, $f_*\Mov (nD)\le\Mov (nf_*D)$.
This proves the statement.
\end{proof}
Example \xref{ex-cremona} shows that in general case inequality
\eqref{eq-push-down} is not an equality.

From Proposition \xref{prop-first-s-decomp} we obtain:
\begin{corollary}[cf. {\cite[Th. 3.5. 3] {Nak}}]
Let $f\colon Y\to X$ be a birational contraction of
$\QQ$-factorial varieties and let $D$ be an effective modulo $\qq$
divisor on $X$. Then
\[
f^*P_{\s} (D)\ge P_{\s} (f^*D).
\]
\end{corollary}

Note that the s-decomposition does not satisfy condition
\eqref{surd}:

\begin{example}
\label{ex-P3-2tochki}
Let $X$ be a blowup of $\PP^3$ in two distinct points $P_1$ and
$P_2$, and let $D$ be the birational transform of a plane passing
through $P_1$ and $P_2$. Since the divisor $D$ is mobile, $N_{\s}
(D) =0$. Now let $f\colon Y\to X$ be the blowup of the birational
transform of the line passing through $P_1$ and $P_2$ and let $E$
be the exceptional divisor. Then $N_{\s} (f^*D) =E\neq f^*N_{\s}
(D)$.
\end{example}

\subsection*{Sectional decomposition}
\begin{definition}[{\cite[\S 2]{Kaw-crep}}]
\label{def-sect-dec}
A decomposition $D=P_{\sect} (D) +N_{\sect} (D)$ is called a
\emph{sectional decomposition} if the following conditions are
satisfied:
\begin{enumerate}
\item
$N_{\sect} (D)\ge 0$;
\item
$P_{\sect} (D)\in\ov\Mv (X)$ (in particular, $P_{\sect}$ is
$\RR$-Cartier);
\item
there is an isomorphism of graded algebras
\begin{equation*}
\label{CKM-ravenstvo}
\sR_XD\simeq\sR_XP_{\sect}.
\end{equation*}
\end{enumerate}
\end{definition}

From definition and Proposition \xref{prop-first-s-decomp} we
immediately obtain:
\begin{proposition}
\label{prop-sect-s}
Let $D$ be an effective modulo $\qq$ divisor on a $\QQ$-factorial
variety. Then an s-decomposition is also a sectional
decomposition. For any sectional decomposition we have
$P_{\sect}\ge P_{\s}$
\end{proposition}

\begin{remarks}
\begin{enumerate}
\item
Note that unlike s- and the $\sigma$-decompositions a sectional
decomposition not necessary unique (in Example \xref{ex-9-points}
there exists infinitely many sectional decompositions $P_{\sect}
(D) =tD$, $0\le t\le 1$, while $P_{\s} (D) =0$).
\item
The set of all sectional of decompositions is closed: if there
exists a sequence $D=P_{\sect}^{(r)}+N_{\sect}^{(r)}$  of
sectional decompositions and the limit
$P=\lim_{r\to\infty}P_{\sect}^{(r)}$ exists, then and $D=P+ (D-P)$
is also a sectional decomposition.
\end{enumerate}
\end{remarks}

\begin{proposition}[{\cite[\S 2] {Kaw-crep}}]
\label{sect-dec-ed}
Let $D$ be a big divisor on a $\QQ$-factorial variety. Then a
sectional decomposition is unique \textup(and
$P_{\sect}=P_{\s}=P_{\sigma}$\textup). Furthermore, by Proposition
\xref{prop-s-decomp-unique} the negative part $N_{\sect} (D)$ is
uniquely defined by the numerical class of $D$.
\end{proposition}
\begin{outline}
According to Proposition \xref{prop-sect-s} it is sufficient to
prove the inequality $P_{\sect}\le P_{\s}$. Note that the divisor
$P_{\sect}$ is big. As in the proof of Proposition
\xref{prop-s-decomp-unique} by the Kodaira lemma we write
$P_{\sect}=H+F$, where $H$ is ample and $F$ is effective. Then for
any $0 <\ep < 1$ the class of $P_{\sect} -\ep F= (1 -\ep)
P_{\sect}+\ep H$ is contained in the cone $\Mv^o (X)$. According
to Example \xref{ex-D=0,D=1}, (ii) we have
\[
P_{\s} (D)\ge P_{\s} (P_{\sect} -\ep F) =P_{\sect} -\ep F.
\]
Passing to the limit, we obtain $P_{\sect}\le P_{\s}$.
\end{outline}

\section{Zariski decomposition in Shokurov's sense and bss-ampleness}
\label{sect-Shok}
\subsection*{Bss-ampleness}
We need the following, almost obvious lemma.

\begin{lemma}
Let $\sD$ be a b-divisor of the fields $K (X)$ and let $L$ be a
divisor on $X$ such that $\sD_X\qr L$. Then there exists a unique
b-divisor $\varLambda (L,\sD)$ such that $\varLambda (L,\sD)_X=L$
and $\varLambda (L,\sD)\qr\sD$.
\end{lemma}

\begin{definition}[{\PLF{Def. 3.2}}, see also {\cite[Prop. 1.10]
{F3}}]
\label{def-bss}
A divisor $D$ on $X$ is said to be \emph{bss-ample} if there
exists a rational $1$-contraction $\alpha\colon X\dasharrow Y$, a
(numerically) ample divisor $H$ on $Y$, and a decomposition
\begin{equation}
\label{eq-diag-bss-d}
D=D^{\m} +E,
\end{equation}
satisfying the following conditions:
\begin{enumerate}
\item
$D^{\m}\qr\alpha^*H$ and there exists a b-semiample b-divisor
$\sH$ such that $\sH_X=\alpha^*H$;
\item
the divisor $E$ is effective and very exceptional with respect to
the map $\alpha$;
\item
for any b-semiample b-divisor $\sL$ such that $\sL_X\le D$ we have
$\sL\le\varLambda (D^{\m},\sH)$.
\end{enumerate}
\end{definition}

\begin{remark}
\begin{enumerate}
\item
According to (iii) the b-divisor $\varLambda (D^{\m},\sH)$ is
uniquely defined (i.e., it depends only on $D$). Therefore
divisors $D^{\m}$, $E$, and the contraction $\alpha$ also are
uniquely defined. The b-divisor $\sH$ is defined up to
$\RR$-linear equivalence.

\item
In the notation of Definition \xref{def-bss} the divisor $D$ is
big if and only if the rational $1$-contraction $\alpha$ is
birational (see \PLF{Prop. 3.20}).
\item
The condition of (iii) in the definition is automatically
satisfied in cases $\dim X -\dim Y\le 1$ and $\dim Y=0$ (see
\PLF{3.4.3}).
\end{enumerate}
\end{remark}

\begin{examples}
\label{ex-2}
{\rm (i)} Any semiample divisor $D$ is bss-ample (we can put $E=0$
and $D^{\m}=D$). A b-semiample divisor $D$ is bss-ample if and
only if there exists the largest b-semiample b-divisor
$\sD^m=\varLambda (D^{\m},\sH)$ with $\sD^m_X=D$.

{\rm (ii)} In the case $\kappa (X, D) =0$ the divisor $D$ is
bss-ample if and only if $D\qr 0$.

{\rm (iii)} Consider the map $\alpha$ from Example
\xref{ex-Pn-bundle} and let $D$ be an $n$-dimensional plane on
$X$. Then $D$ is a bss-ample divisor because
$D=\alpha^*\OOO_{\PP^1} (1)$.

{\rm (iv)} A divisor $D$ with $\kappa (X, D) > 0$ on a surface is
bss-ample if and only if the positive part $P (D)$ of the
classical Zariski decomposition is semiample. In this instance
decomposition \eqref{eq-diag-bss-d} coincides with the Zariski
decomposition, i.e., $N (D) =E$.

Indeed, by Zariski's Main Theorem the map $\alpha$ is a morphism.
In the case $\dim Y=2$ the contraction $\alpha$ is birational and
the decomposition $D=D^{\m} +E$ coincides with the Zariski
decomposition by definition and because the intersection matrix of
exceptional divisors is negative definite. If $\dim Y=1$, then the
intersection matrix $(E_i\cdot E_j)$ is negative definite by the
property of very exceptionality and semi-negativity of the
intersection matrix in fibers of $X\to Y$.

Conversely, if $P (D)$ is a semiample divisor, then it defines a
contraction $\alpha\colon X\to Y$, for which the divisor $N (D)$
is exceptional. Since the intersection matrix of components of $N
(D)$ is negative definite, in case $\dim Y=1$, a fiber of the
contraction $\alpha$ cannot be contained in $\Supp (N (D))$.
\end{examples}

\begin{remark}
\label{rem-bss}
Let $D$ be a bss-ample divisor on $X$ and let $\alpha\colon X\bir
Y$ be the corresponding rational $1$-contraction. Then there
exists a diagram
\begin{equation}
\label{diag-bss-bss}
\begin{diagram}
&&W&&
\\
&\ldTo^{\scriptstyle{g}}&&\rdTo^{\scriptstyle{h}}&
\\
X&&\rDashto^{\scriptstyle{\alpha}}&&Y
\\
\end{diagram}
\end{equation}
and a numerically ample divisor $H$ on $Y$ such that
\[
g\1D+F_W\qr h^*H+E_W,
\]
where divisors $F_W$ and $E_W$ are exceptional on $X$ and $Y$,
respectively, $E_W$ is effective, and the divisor $g (E_W)$ is
very exceptional on $Y$. It is clear that $\sH=\overline{h^*H}$.

Note however that divisors $F_W$ and $E_W$ are defined
ambiguously. We can take as $E_W$ the birational transform of $E$
and in this instance $E_W$ will be very exceptional on $Y$.
\end{remark}

Obviously any bss-ample divisor is effective (more precisely,
$D\qr D'$, where $D'\ge 0$). However, the converse is not always
true: the divisor $D$ from Example \xref{ex-ell-Zar} is effective
and nef but it is not bss-ample. Indeed, this divisor is not
semiample (because its restriction to $\widetilde C$ is a
numerically trivial divisor which is not a torsion).

\begin{proposition}
\label{prop-bss-pairs}
Let $D$ be a big divisor. Assume that the $D$-MMP holds
\textup(including $D$-abundance conjecture\textup). Then $D$ is
bss-ample.
\end{proposition}
The conditions of the proposition are satisfied when $\dim X\le 3$
(and in any dimension modulo LMMP) in most important cases (see
\xref{p-exist} below).

\begin{proof}
First replace $(X, D)$ with its $\QQ$-factorialization, and then
run the $D$-MMP. It is clear that bss-ampleness is invariant under
birational isomorphisms in codimension $1$. Thus it is sufficient
to show that bss-ampleness is preserved under divisorial
contractions:

\begin{lemma}
Let $\var\colon X\to X'$ be a divisorial extremal $D$-negative
contraction and let $D'=\var_*D$. Assume that $D'$ is big and
bss-ample. Then so is $D$.
\end{lemma}
\begin{proof}
Let $S$ be an (irreducible) $\var$-exceptional divisor. Then
$D=\var^*D'+aS$, where $a\ge 0$ according to Lemma \xref{11}.
Consider Hironaka's hut
\begin{equation*}
\begin{diagram}
&&W&&
\\
&\ldTo^{\scriptstyle{\delta}}&
\dTo^{\scriptstyle{\pi}}&\rdTo^{\scriptstyle{\tau}}&
\\
X&\rTo^{\scriptstyle{\var}}&X'& \rDashto^{\scriptstyle{\beta}}&Y
\\
\end{diagram}
\end{equation*}
and let $\alpha\eqdef\beta\comp\var$. Then $\alpha$ is a rational
$1$-contraction. Furthermore, since $D'$ is big, the map $\alpha$
is birational. We have $D'=D^{\prime\m}+E'$, where
$D^{\prime\m}\qr\beta^*H$ and the divisor $H$ is ample. Put
\[
E=\var^*\pi_*\tau^*H -\delta_*\tau^*H+\var^*E' +aS,\qquad
D^{\m}=D-E.
\]
Then
\begin{multline*}
D^{\m}=\var^*D'+aS - (\var^*\pi_* \tau^*H -\delta_*\tau^*H+
\var^*E' +aS)=
\\
\var^* (D' -\beta*H'-E')+\alpha^*H\qr\alpha^*H.
\end{multline*}
Further,
\[
\var_*E=\pi_*\tau^*H -\var_*\delta_*\tau^*H+E'=E'\ge 0.
\]
One can also see that the divisor $E$ is exceptional on $Y$.
Finally, by Lemma \xref{11} the divisor $\var^*\pi_*
\tau^*H+\var^*E' +aS -\delta_*\tau^*H$ is effective.
\end{proof}

After a finite number of contractions and flips $X\bir
X_1\bir\cdots\bir X_n$ we obtain a model $(X_n, D_n)$, on which
the divisor $D_n$ is nef (and big). By our abundance hypothesis
the divisor $D_n$ is semiample. Hence it is bss-ample.
\end{proof}

Note that it is not possible to omit the bigness condition in
Proposition \xref{prop-bss-pairs}:

\begin{example}
\label{ex-bss-Sh-not}
Let $X=\PP^1\times\FF_1$, let $F$ be a fiber of the projection
$\PP^1\times\FF_1\to\PP^1$, and let $G=\PP^1\times C_0\subset X$,
where $C_0$ is the negative section of $\FF_1$. Consider the
divisor $D=F+G$. Then the contraction $\var\colon
X\to\PP^1\times\PP^2$ is a (unique) step of $D$-MMP and the
divisor $\var_*D$ is nef. However, $D$ is not bss-ample. Indeed,
suppose that there is a decomposition such as in
\eqref{eq-diag-bss-d}: $D=D^{\m}+E$. Since the divisor $D^{\m}$ is
b-semiample, we have that there exist at most a finite number of
curves negatively intersecting $D^{\m}$. On the other hand, $X$ is
a smooth Fano variety. Hence the divisor $D^{\m}$ is nef and the
rational $1$-contraction $\alpha$ is a morphism. Further, $X$ is a
quasihomogeneous variety. We obtain that for $E$ there is only one
possibility $E=\lambda G$, $\lambda\ge 0$. Then $D^{\m}=F+\lambda
G$, where $0\le\lambda\le 1$. However, $G$ is a fixed component of
the linear system $\lins{n (F+\lambda G)}$ for $\lambda > 0$.
Therefore, $\lambda=0$, $D^{\m}=F$, and the rational
$1$-contraction $\alpha$ coincides with the projection
$X\to\PP^1$. But the divisor $G$ is not exceptional on $\PP^1$, a
contradiction.

Note that $(X,\Delta)$ is $0$-pair for a suitable boundary
$\Delta$.
\end{example}

\subsection*{Zariski decomposition in Shokurov's sense}

\begin{definition}
\label{def-Zar_Shok}
A decomposition $D=D^{\m}+D^{\e}$ is called a \emph{Zariski
decomposition in sense of Shokurov} (or simply \emph{Shokurov
decomposition}) if the following conditions are satisfied:
\begin{enumerate}
\item
$D^{\e}\ge 0$;
\item
there exists a b-semiample b-divisor $\sD^{\m}$ such that
$\sD^{\m}_X=D^{\m}$;
\item
for any b-semiample b-divisor $\sL$ such that $\sL_X\le D$ we have
$\sL\le\sD^{\m}$.
\end{enumerate}
\end{definition}

Thus any bss-ample divisor has a Shokurov decomposition (and it
coincides with \eqref{eq-diag-bss-d}). In general case, the
converse not always true:

\begin{example}
Let $X$ and $D$ be such as in Example \xref{ex-bss-Sh-not}. Then
we have the following Shokurov decomposition:
\[
D^{\m}=P_{\s} (D) =F,\qquad D^{\e}=G.
\]
However, the divisor $D$ is not exceptional for the morphism
$X\to\PP^1$ defined by the linear system $\lins{F}$. Therefore $D$
is not bss-ample.
\end{example}

Immediately from definition we obtain the following properties:

\begin{pusto}
\label{pusto-propert-Sh}
Let $D$ be an effective modulo $\QQ$-linear equivalence divisor.
\begin{enumerate}
\item
If a Shokurov decomposition of $D$ exists, then $D^{\m}\ge P_{\s}
(D)$. If furthermore $D^{\m}$ is a $\QQ$-divisor, then
$D^{\m}=P_{\s} (D)$ (see Proposition \xref{prop-s-dec-semiample}).
\item
If a Shokurov decomposition for of $D$ exists, then
$\sR_XD^{\m}=\sR_XD=\sR_XP_{\s} (D)$ (see Proposition
\xref{prop-first-s-decomp}).
\end{enumerate}
\end{pusto}

\begin{remark}
\label{remark-propert-Sh}
Let $P_{\s} (D)$ be a b-semiample $\QQ$-divisor. Then for the
decomposition $D=P_{\s}+N_{\s}$ conditions (i) and (ii) of
Definition \xref{def-Zar_Shok} are satisfied while (iii) is not
necessary satisfied. Instead the following weaker condition holds:
\begin{enumerate}
\item[(iii)${}'$]
for any b-semiample divisor $L$ such that $L\le D$ we have $L\le
D^{\m}$.
\end{enumerate}
Indeed, it is sufficient to show that for any b-semiample divisor
$L$ such that $L\le D$ the inequality $P_{\s}\ge L$ holds.

Assume that $L\nleq P_{\s}$. Write $P_{\s}=\sum\alpha_i L_i$ and
$L=\sum\beta_i L_i$, where the $L_i$ are integral b-free divisors,
$\alpha_i\in\QQ_{> 0}$, and $\beta_i\in\RR_{> 0}$. Consider the
finite dimensional real vector space $V$ generated by the
components of divisors $P_{\s}$ and $L$. The conditions $\Min
(P_{\s}, L)\le x\le\Max (P_{\s}, L)$ define a parallelepiped
$\mathfrak K$ (possible not maximal dimension) in $V$ and the
b-semiample divisors $L_i$ generate a convex rational polyhedron
$\mathfrak R\subset V$ containing the diagonal $[P_{\s}, L]$ of
$\mathfrak K$. So it is easy to see that there exists a rational
divisor $L'$ into the interior of the set $\mathfrak
K\cap\mathfrak R$. From Proposition \xref{prop-s-dec-semiample} we
obtain $L'\le P_{\s}$. Contradiction with the fact that $L\nleq
P_{\s}$.
\end{remark}

From Theorem \xref{th-Lim-Crit} and \xref{remark-propert-Sh} we
obtain:

\begin{corollary}[cf. \PLF{Remark 3.30, Example 4.30}]
\label{corollary-last}
Suppose that a divisorial algebra $\sR_XD$ nontrivial and finitely
generated. Then there exists a Shokurov decomposition \textup
(furthermore, $D^{\m}=P_{\s}$ is a $\QQ$-divisor\textup).
\end{corollary}

Remark \xref{remark-propert-Sh} and Corollary
\xref{corollary-last} are particular cases of \PLF{Example 4.30}
(see also Proposition \xref{proposition-b-divisors}).

\begin{proposition}
\label{prop-Sh-equival-1}
Let $D$ be a big divisor on a $\QQ$-factorial variety $X$. If a
Shokurov decomposition for $D$ exists, then $D^{\m}=P_{\s}$. In
particular, $D^{\e}$ depends only on the class of $D$ modulo
numerical equivalence.
\end{proposition}

\begin{proof}
Suppose that $D^{\m}\gneqq P_{\s} (D)$. Similar to the proof of
Proposition \xref{sect-dec-ed}, we can write $D^{\m}=H+F$, where
$H$ is ample and $F$ is an effective divisor. Then for any $0 <\ep
< 1$ the class of the divisor $D^{\m} -\ep F= (1 -\ep) D^{\m}+\ep
H$ is contained in $\Mv^o (X)$. According to Example
\xref{ex-D=0,D=1}, (ii) we have
\[
P_{\s} (D)\ge P_{\s} (D^{\m} -\ep F) =D^{\m} -\ep F.
\]
Passing to the limit, we obtain $D^{\m}\le P_{\s}$, a
contradiction.
\end{proof}

\begin{example}
The integral divisor $D$ from Example \xref{ex-ell-Zar} has no
Shokurov decompositions. In fact, we can assert that if a divisor
$D$ on variety of any dimension is nef and for any $n\in\NN$ the
linear system $\lins{nD}$ has fixed components of bounded
multiplicity, then $D$ has no Shokurov decompositions.
\end{example}

\begin{example}
Let $D$ be a divisor from Example \xref{ex-9-points}. Since
$D^2=0$ and $\dim\lins{nD}=0$ for all $n\in\NN$, we have that $D$
is not semiample. Therefore $D$ has the following Shokurov
decomposition $D^{\m}=0$ and $D^{\e}=D$. It does not coincide with
the classical Zariski decomposition $N (D) =0$.
\end{example}

In the two-dimensional case there is a criterion of the existence
of Shokurov decomposition (see Example \xref{ex-2} and Proposition
\xref{prop-Zar-surf}):
\begin{proposition}
\label{prop-surf-Zar-Sh}
Let $X$ be a surface with $\QQ$-factorial singularities and let
$D$ be a divisor on $X$ with the \textup(classical\textup) Zariski
decomposition $D=P+N$. Then the following conditions are
equivalent\textup:
\begin{enumerate}
\item
the divisor $P$ is semiample;
\item
the divisor $D$ is bss-ample;
\item
there exists a Shokurov decomposition $D=D^{\m}+D^{\e}$ such that
$D^{\e}=N$ and $\sD^{\e}=\overline N$.
\end{enumerate}
\end{proposition}

\begin{example}
If in the conditions of Proposition \xref{prop-surf-Zar-Sh} on $X$
there exists a Shokurov decomposition $D=D^{\m}+D^{\e}$ and
$\kappa (X, D)\ge 1$, then $D^{\e}= N (D)$.

Indeed, by Proposition \xref{prop-Zar-surf} we always have $P
(D)\ge D^{\m}\ge P_{\s}$ and according to Remark \xref{rem-s-Zar}
the equality holds.
\end{example}

\begin{theorem}[\PLF{Th. 3.33}, cf. Proposition {\xref{prop-bss-pairs}},
Theorem {\xref{th-CKM-logcanonical}}]
\label{th-LMMP-Shok}
Let $(X, B)$ be a $0$-pair, i.e., a Kawamata log terminal pair
with $K_X+B\equiv 0$ and let $D$ be an effective modulo
$\RR$-linear equivalence divisor. Assume the LMMP
\textup(including the abundance conjecture\textup). Then there
exists a Shokurov decomposition for the divisor $D$.
\end{theorem}
In the conditions of this theorem there exist (and coincide
between each other) all other Zariski decompositions (see below)
for the divisor $D$.

\section{Zariski decomposition in Fujita's sense}
\label{sect-Fudj}
Fujita noticed that property \xref{prop-Zar-surf} characterizes
Zariski decomposition:

\begin{definition}
\label{def-Fujita-new}
A decomposition $D=P_{\F}+N_{\F}$ is called a \emph{Zariski
decomposition in Fujita sense} (or simply \emph{Fujita
decomposition}) if
\begin{enumerate}
\item
$N_{\F}\ge 0$;
\item
$P_{\F}$ is nef;
\item
for any b-nef divisor $L$ such that $L\le D$ we have $L\le
P_{\F}$.
\end{enumerate}
\end{definition}

It is clear that we may set a question about the existence of such
a decomposition only for pseudo-effective divisors.

\begin{remarks}
\begin{enumerate}
\item
It follows immediate from the definition that the negative part of
a Fujita decomposition depends only on the class of numerical
equivalence of $D$.
\item
It is clear that $P_{\F}\ge M_n/n$. Therefore, $P_{\F}\ge P_{\s}$.
So,
\begin{equation}
\label{eq-f-s}
\sR_XD=\sR_XP_{\F}=\sR_XP_{\s}.
\end{equation}
Example \xref{ex-9-points} shows that the equality $P_{\F}=
P_{\s}$ is not always true (even if the divisor $P_{\s}$ is nef).
Nevertheless $D=P_{\F}+N_{\F}$ is a $\sigma$-decomposition for any
pseudo-effective divisor $D$ (see \cite{Nak}).
\item
In dimension $2$ the Fujita decomposition coincides with the
classical Zariski decomposition (see Proposition
\xref{prop-Zar-surf}).
\end{enumerate}
\end{remarks}

\begin{example}
Let $D$ be a big bss-ample $\RR$-Cartier divisor. We make use the
notation of Remark \xref{rem-bss}. Then in diagram
\eqref{diag-bss-bss} there exists a unique decomposition
\begin{equation}
\label{eq-R-Cartier}
g^*D =P_{\F}+N_{\F},
\end{equation}
where $P_{\F}\qr h^*H$ and the divisor $N_{\F}$ is effective and
(very) exceptional on $Y$. This decomposition is a Fujita
decomposition (with semiample positive part).

Indeed, in this instance $h$ is a birational contraction. Put
$N_{\F}= E_W-F_W$. Then the divisor $-N_{\F}\qr h^*H-g^*D$ is nef
over $X$. By Lemma \xref{11} we have $N_{\F}\ge 0$. It is clear
that $N_{\F}$ is exceptional on $Y$. Assume that a divisor $L\le
D$ is nef. Then the divisor $L-P_{\F}\qr L - h^*H$ is nef over
$Y$. From Lemma~ \xref{11} we obtain that $L\le P_{\F}$. Thus
$g^*D\qr P_{\F}+N_{\F}$ is a Zariski decomposition in Fujita
sense. It is unique.

The statement is no longer true if we omit the condition for $D$
to be big (see Example \xref{ex-Pn-bundle})
\end{example}

\begin{proposition}
\label{prop-def-Fujita-new}
Let $D$ be an effective $\RR$-Cartier divisor. Assume that the
decomposition $D=P_{\F}+N_{\F}$ satisfies conditions {\rm
(i)-(ii)} of Definition \xref{def-Fujita-new}. Then {\rm (iii)} is
equivalent to the following
\begin{enumerate}
\item[(iii)${}'$]
for any birational contraction $f\colon\hat X\to X$ and for any
nef divisor $F$ on $\hat X$ such that $F\le f^*D$ we have $F\le
f^*P_{\F}$.
\end{enumerate}
\end{proposition}
\begin{proof}
\textbf{(iii) $\Rightarrow$ (iii)${}'$.} Since $-(f^*P_{\F}-F)$ is
nef over $X$ and $f_* (f^*P_{\F}-F)\ge 0$, we have $f^*P_{\F}-F\ge
0$. (We used that $f_*F$ is b-nef and $f_*F\le D$, so $f_*F\le
P_{\F}$).

\textbf{(iii) ${}'$ $\Rightarrow$ (iii).} Since $L$ is b-nef,
there exists a birational contraction $f\colon\hat X\to X$ and a
nef divisor $\hat L$ on $\hat X$ such that $f_*\hat L=L$. Further,
$\hat L-f^*D$ is nef over $X$ and $f_* (f^*D -\hat L)\ge 0$.
Hence, $f^*D -\hat L\ge 0$ and $L\le P_{\F}$.
\end{proof}

However, a Fujita decomposition does not always exists that shows
the following simple example:

\begin{example}
\label{ex-blow-up-P3}
Consider the linear system $\lins{D}$ on a nonsingular projective
variety $X$ such that
\begin{enumerate}
\item
$\lins{D}$ has no fixed components;
\item
$D$ is not nef.
\end{enumerate}
For example, similar to Example \xref{ex-P3-2tochki} we may take
as $X$ a blowup of $\PP^3$ in two distinct points $P_1$ and $P_2$,
and as $D$ the birational transform of a plane passing through
$P_1$ and $P_2$. Assume that there exists a decomposition $N_{\F}
(D)$ in Fujita's sense. Then $P_{\F}\ge M_1=D$, i.e., $P_{\F}=D$,
a contradiction.
\end{example}

\begin{pusto}
\label{p-exist}
The above shows that it is more naturally to construct a Fujita
decomposition for a pull-back $f^*D$ of $D$ under some birational
contraction $f\colon Y\to X$. However even in such a stating, the
problem of the existence of a Fujita decomposition fails (see
\cite{Nak}). Nevertheless, its positive decision is expected (for
a pseudo-effective divisor $D$) in the most important cases:
\begin{enumerate}
\item
if $(X, B)$ is $0$-pair (see Theorem \xref{th-LMMP-Shok});
\item
if $(X, B)$ is a Fano log variety;
\item
if $D=K_X+B$, where the pair $(X, B)$ has log canonical
singularities and $\kappa (X, K_X+B)\ge 0$ (cf. Theorem
{\xref{th-CKM-logcanonical}}).
\end{enumerate}
Argument justifying this hope is, for example, fact proved by
Kawamata \cite[Prop. 5] {Ka}: for any effective Cartier divisor
$D$ on a toric variety $X$ there exists a toric birational
contraction $f\colon Y\to X$ such that the divisor $f^*D$ has a
Fujita decomposition and this decomposition coincides with a
Shokurov decomposition and with an CKM-decomposition (see below).
\end{pusto}

\begin{remark}
\label{zam-Fudj-b-div}
It is clear that the Fujita decomposition satisfies condition
\eqref{surd}. Therefore (see Remark \xref{zam-b-div}) if it exists
on $Y$ for the divisor $\sD_Y$, where $\sD\in\bCDiv_{\RR} (K
(Y))$, then there exists a b-divisor $\sN (\sD)\in\bCDiv_{\RR} (K
(Y))$ such that $\sN (\sD)_{Y'}$ is a decomposition in Fujita
sense for any model $Y'$ of the field $K (Y)$ dominating $Y$. In
other words a Fujita decomposition is unique \emph{in the
birational sense}.
\end{remark}

\begin{definition}
\label{def-Fujita-ger}
A decomposition $D=P_{\gF}+N_{\gF}$ is called \emph{generalized
Fujita decomposition} if
\begin{enumerate}
\item
$N_{\gF}\ge 0$;
\item
$P_{\gF}$ is b-nef;
\item
for any b-nef divisor $L$ such that $L\le D$ we have $L\le
P_{\gF}$.
\end{enumerate}
\end{definition}

Another way to generalize Fujita decomposition is to consider the
straightforward analog of Definition \xref{def-Zar_Shok}:
necessary to claim the existence of a b-nef b-divisor $\sP$ such
that $\sP_X=P$ and satisfying condition \xref{def-Zar_Shok} (iii).

\begin{proposition}
\label{prop-def-Fujita-gen}
Let $D$ be an effective $\RR$-Cartier divisor. Assume that there
exists a birational contraction $f\colon\hat X\to X$ such that the
divisor $D^*=f^*D$ has a Fujita decomposition $D^*= P_{\F}
(D^*)+N_{\F} (D^*)$. Then $D=f_*P_{\F} (D^*)+f_* N_{\F} (D^*)$ is
a generalized Fujita decomposition.
\end{proposition}
\begin{proof}
By construction $f_*N_{\F} (D^*)$ is effective and $f_*P_{\F}
(D^*)$ is b-nef. Let $L$ be a b-nef divisor on $X$ such that $L\le
D$. We may assume that there exists a nef divisor $\hat L$ on
$\hat X$ such that $f_*\hat L=L$. Then $-(f^*D -\hat L)$ is nef
over $X$ and $f_* (f^*D -\hat L)\ge 0$. By Lemma \xref{11} we have
$f^*D\ge\hat L$. But then $P_{\F} (D^*)\ge\hat L$ and $f_*P_{\F}
(D^*)\ge L$.
\end{proof}

\begin{proposition}
\label{prop-Fujita-gen-1}
Let $D$ be an effective divisor. Assume that there exists a
generalized Fujita decomposition $D=P_{\gF}+F_{\gF}$. Then
\begin{enumerate}
\item
$\sR_XD=\sR_XP_{\gF}$;
\item
$P_{\gF}\ge P_{\s} (D)$. If furthermore the divisor $D$ is big and
the variety $X$ is $\QQ$-factorial, then $P_{\gF}=P_{\s} (D)$.
\end{enumerate}
Conversely, if the divisor $D$ is big, $X$ is $\QQ$-factorial, and
$P_{\s} (D)$ is b-nef, then there exists the generalized Fujita
decomposition for $D$ and $P_{\gF}=P_{\s} (D)$.
\end{proposition}
\begin{proof}
The proof of (i) is similar to the proof of Proposition
\xref{prop-first-s-decomp}. We prove (ii). Since $M_n\le nP_{\gF}$
for all $n\in\NN$, we have $P_{\s} (D)\le P_{\gF}$. Assume now
that the divisor $D$ is big. By (i) so is $P_{\gF}$. Fix $\ep> 0$.
By the Kodaira lemma $P_{\gF}=A+F$, where $A$ is ample and $F$ is
an effective divisor. It is clear that $P_{\gF} -\ep F= (1 -\ep)
P_{\gF}+\ep A$ is a b-semiample divisor and its class is contained
in the open cone $\Mv^o (X)$. Similar to the proof of Proposition
\xref{prop-Sh-equival-1} we have
\[
P_{\s} (D)\ge P_{\s} (P_{\gF} -\ep F) =P_{\gF} -\ep F.
\]
Passing to the limit, we obtain $P_{\gF}\le P_{\s} (D)$.

Finally, let $L\le D$ be a b-nef divisor. We prove that $L\le
P_{\s} (D)$. By the Kodaira lemma $D=A+F$, where the divisor $A$
is ample and $F$ is effective. Take a sufficiently small rational
$\ep > 0$. Then $(1+\ep) D\ge L+\ep A$. As above
\[
(1+\ep) P_{\s} (D) =P_{\s} ((1+\ep) D)\ge P_{\s} (L+\ep A) =L+\ep
A.
\]
Hence, $L\le P_{\s} (D)$. This proves the proposition.
\end{proof}

\pusto{}\ \textbf{Zariski decomposition for varieties of
intermediate Kodaira dimension.}

\begin{proposition}[{\cite[Prop. 1.24, 1.10] {F3}}]
Let $f\colon X\to Z$ be a contraction of nonsingular varieties,
let $D$ be a divisor on $Z$, and let $R$ be a very exceptional
divisor on $X$. Then the pull-back of $f^*D+R$ under some
birational contraction $g\colon X'\to X$ has a Zariski
decomposition in Fujita sense if and only if this decomposition
exists for the pull-back of $D$ under some birational contraction
$h\colon Z'\to Z$. Furthermore, the divisor $P_{\F} (g^*
(f^*D+R))$ is the pull-back of $P_{\F} (h^*D)$ for a suitable
choice of $g$ and $h$: if the diagram
\[
\begin{CD}
X@ < g << X'
\\
@VfVV @V{f'}VV
\\
Z@ < h << Z'
\end{CD}
\]
is commutative, then $P_{\F} (g^* (f^*D+R)) =f^{\prime *}P_{\F}
(h^*D)$.
\end{proposition}

Using this statement and a formula for canonical divisor of
elliptic fibrations Fujita proved the following theorem.

\begin{theorem}[{\cite{F3}}]
Let $f\colon X\to Z$ be a contraction of a three-dimensional
variety with general fiber being an elliptic curve. Assume that
$\kappa (X, K_X)\ge 0$. Then there exists a birational contraction
$g\colon X'\to X$ such that $g^*K_X$ has a Zariski decomposition
in Fujita sense $P_{\F}+N_{\F}$ with semiample positive part
$P_{\F}$ \textup(and it coincides with the Shokurov
decomposition\textup). In particular, the canonical algebra $\sR
K_X$ is finitely generated.
\end{theorem}

Later the last fact was generalized in \cite{FM}:

\begin{theorem}[{\cite{FM}}]
Let $(X, B)$ be a Kawamata log terminal pair with $\kappa (X,
K_X+B) =l\ge 0$, where $B$ is a $\QQ$-boundary. Then there exists
an $l$-dimensional Kawamata log terminal pair $(Z,\Delta)$ with
$\kappa (Z, K_Z+\Delta) =l$ such that log canonical algebras $\sR
(K_X+B)$ and $\sR (K_Z+\Delta)$ are quasi-isomorphic \PLF{Def.
4.3}. In particular, questions about finite generation of these
algebras are equivalent \PLF{Th. 4.6}.
\end{theorem}

\begin{corollary}
Let $(X, B)$ be a Kawamata log terminal pair with $0\le\kappa (X,
K_X+B)\le 3$. Then the log canonical algebra $\sR (K_X+B)$ is
finitely generated.
\end{corollary}

\pusto{}\ \textbf{An example of a divisor with non-rational
Zariski decomposition.}
\label{ex-nonrat}
Following \cite{Cu}, we present an example of a divisor with a
Zariski decomposition in Fujita sense having irrational
coefficients.

Let $E$ be an elliptic curve with $\mt{End} (E)\simeq\ZZ$ and let
$S=E\times E$. Pick a point $P\in E$ and put $E_1=\{P\}\times E$
and $E_2=E\times\{P\}$. Then the diagonal $\Delta$ and the curves
$E_1$ and $E_2$ form the basis of the space $\NS_{\RR} (S)$. Put
$\Delta'=\Delta-E_1-E_2$. Then $\Delta', E_1, E_2$ is also a
basis. In this basis the quadratic form $x^2$ has the following
form
\[
(\alpha_1 E_1+\alpha_2 E_2+\beta\Delta')^2= 2\alpha_1 \alpha_2
-2\beta^2.
\]
Therefore the cone of ample divisors in $\NS_{\RR} (S)$ is defined
by conditions
\begin{equation}
\label{eq-alpha}\alpha_1\alpha_2
-\beta^2 > 0,\qquad\alpha_1,\alpha_2 > 0.
\end{equation}

Consider the $\PP^1$-bundle $\pi\colon X=\PP (\Oh_S
(\Delta')\oplus\Oh_S)\to S$. Identify $S$ with the zero section.
Then $\Oh_S (S)=\Oh_S (\Delta')$. Put $S_1=\pi^*E_1$ and
$S_2=\pi^*E_2$. Let $H$ be an integral ample divisor on $X$ and
let $L\eqdef H|_S$.

Consider the divisor $G (\alpha_1,\alpha_2)\eqdef L+\alpha_1
E_1+\alpha_2 E_2$. For positive $\alpha_1,\alpha_2$, define the
function
\[
\gamma (\alpha_1,\alpha_2)\eqdef\sup\{\beta\mid G
(\alpha_1,\alpha_2)+\beta\Delta'\quad\text{is nef}\}
\]
Using \eqref{eq-alpha}, choose $\alpha_1,\alpha_2\in\QQ$,
$\alpha_1,\alpha_2 > 0$ so that $\gamma\eqdef\gamma
(\alpha_1,\alpha_2)$ is irrational. If $\delta >\gamma$, then $G
(\alpha_1,\alpha_2)+\delta\Delta'$ is not nef. Therefore there
exists an irreducible curve $\Gamma$ having negative intersection
with $G (\alpha_1,\alpha_2)+\delta\Delta'$. Since $S$ is an
abelian surface, there exists a family of such curves
$\{\Gamma_{\lambda}\}$ on $S$:
\begin{equation}
\label{eq-family}
(G (\alpha_1,\alpha_2)+ \gamma\Delta')\cdot\Gamma_{\lambda} < 0.
\end{equation}
On the other hand, if $0 <\delta <\gamma$, then $G
(\alpha_1,\alpha_2)+ \delta\Delta'$ is ample.

Put $D\eqdef B+rS$, where $B\eqdef H+\alpha_1S_1+ \alpha_2S_2$ and
take $r\in\QQ$ so that $r >\gamma$.

\begin{claim}
Let $D$ be the \textup(effective and big\textup) divisor
constructed above. Then $N_{\F} (D)=(r -\gamma) S$ is a Fujita
decomposition. In particular, for any birational contraction
$f\colon Y\to X$ the pull-back $f^*D$ has no Fujita decompositions
with rational coefficients.
\end{claim}
\begin{proof}
First we show that the divisor $P_{\F}= D-N_{\F} (D)$ is nef.
Indeed, suppose that $(B+\gamma S)\cdot C < 0$ for some
irreducible curve $C$. Since $B$ is ample, $S\cdot C < 0$.
Therefore, $(G (\alpha_1,\alpha_2)+\gamma\Delta')\cdot C < 0$, a
contradiction with our choice of $\gamma$.

Further, let $f\colon Y\to X$ be a birational contraction and let
$F$ be an effective divisor on $Y$ such that $f^*D-F=f^*
(P_{\F}+N_{\F})-F$ is nef. From \eqref{eq-family} we have $(G
(\alpha_1,\alpha_2)+ (\gamma+\ep)\Delta') \cdot \Gamma_{\lambda} <
0$ for any $\ep > 0$. Hence, $(P_{\F}+\ep S)\cdot \Gamma_{\lambda}
< 0$. Let $\{\Gamma_{\lambda}'\}$ be a family of irreducible
curves on $Y$ dominating $\{\Gamma_{\lambda}\}$. Then $f^*
(P_{\F}+\ep S)\cdot\Gamma_{\lambda}' < 0$. Since $B$ is ample,
$f^*S\cdot\Gamma_{\lambda}' < 0$. Therefore the family
$\{\Gamma_{\lambda}'\}$ cover the birational transform $\widetilde
S$ of $S$. On the other hand,
\[
0 < (f^*D-F)\cdot\Gamma_{\lambda}'= (f^* (P_{\F}+\ep S+ (r -\gamma
-\ep) S)-F)\cdot\Gamma_{\lambda}'.
\]
Hence,
\[
-(r -\gamma -\ep) S\cdot\Gamma_{\lambda}' < f^* (P_{\F}+\ep
S-F)\cdot\Gamma_{\lambda}' <-F\cdot\Gamma_{\lambda}'.
\]
This means that $f_*F > (r -\gamma -\ep) S$. Since $\ep$ is an
arbitrary positive, $f_*F\ge (r -\gamma) S$. Finally, the
inequality $F\ge f^* (r -\gamma) S$ follows from $f_*F\ge (r
-\gamma) S$ and the fact that the divisor $f^* (r -\gamma) S-F$ is
nef over $X$ (see Lemma~ \xref{11}).
\end{proof}

Nakayama \cite{Nak} using similar construction, constructed an
example of a big (integral) divisor on a nonsingular variety such
that its pull-back under any blowup has no Fujita decompositions
(as well as CKM decompositions, see below).

\section{Zariski decomposition in CKM's sense}

\begin{definition}[\cite{Cu}, \cite{Ka}, \cite{Mo1}]
A decomposition $D=P_{\CKM} (D) +N_{\CKM} (D)$ is called a
\emph{Zariski decomposition in sense Cutkosky-Kawamata-Moriwaki}
(or simply \emph{CKM-decomposition}) if the following conditions
are satisfied:
\begin{enumerate}
\item
$N_{\CKM} (D)\ge 0$;
\item
the divisor $P_{\CKM} (D)$ is nef;
\item
there is an isomorphism of graded algebras
\begin{equation*}
\sR_XD\simeq\sR_XP_{\CKM} (D).
\end{equation*}
\end{enumerate}
\end{definition}

\begin{remark}
Since there is the following embedding of cones
\[
\mt{Nef} (X)=\ov{\mt{Amp}} (X)\subset\ov{\Mv} (X),
\]
we have that an CKM-decomposition is also a sectional
decomposition. If $D=P_{\sect} (D) +N_{\sect} (D)$ is a sectional
decomposition (such as in \xref{def-sect-dec}), then it is an
CKM-decomposition if and only if the divisor $P_{\sect} (D)$ is
nef.
\end{remark}

\begin{remark}
\label{l-6. 4}
If $D$ is $\RR$-Cartier, then a decomposition $N_{\CKM} (D)$
satisfies condition \eqref{surd}.

Indeed, let $f\colon Y\to X$ be a birational contraction. It is
sufficient to show the existence of isomorphisms
\[
\begin{array}{c}
H^0 (\Oh_X (nP_{\CKM}))\simeq H^0 (\Oh_Y (f^*nP_{\CKM})),
\\
H^0 (\Oh_X (nD))\simeq H^0 (\Oh_Y (f^*nD)),
\end{array}
\]
that follows from the fact that $f$ is a contraction.
\end{remark}

\begin{proposition}
\label{F->CKM}
\begin{enumerate}
\item
Let $D$ be a big divisor on a $\QQ$-factorial variety and let
$D=P_{\F} (D) +N_{\F} (D)$ be its Fujita decomposition. Then it is
also an CKM-decomposition.
\item
Let $D$ be a big $\RR$-Cartier divisor. If a decomposition
$D=P_{\CKM} (D) +N_{\CKM} (D)$ exists, then it is a Fujita
decomposition.
\end{enumerate}
\end{proposition}
\begin{proof}
The statement of (i) follows by \eqref{eq-f-s} and (ii) follows by
Propositions \xref{sect-dec-ed} and \xref{prop-Fujita-gen-1}.
\end{proof}

\begin{remark}
\label{rem-CKM-ne-ed}
In general case (i.e., if $D$ is not big) the CKM-decomposition is
not unique and does not coincide with the Fujita decomposition
(and even with the classical Zariski decomposition). For example,
the divisor $D$ from Example \xref{ex-9-points} has infinitely
many CKM-decompositions: $N_{\CKM} (D) =tD$, $0\le t\le 1$.
\end{remark}

\begin{example}
An CKM-decomposition on a surface with $\QQ$-factorial
singularities satisfies condition (iii) of Theorem \xref{main} and
instead of (ii) we have only
\begin{enumerate}
\item[(ii)${}'$]
the matrix $(N_i\cdot N_j)$ is seminegative definite.
\end{enumerate}
Indeed, we may assume that $P_{\CKM}^2=0$. If $P_{\CKM}\cdot N_i
> 0$, then $(P_{\CKM}+ \ep N_i)^2= \ep (2P_{\CKM}\cdot N_i+\ep
N_i^2) > 0$ for some $0 <\ep\ll 1$, i.e., the divisor
$P_{\CKM}+\ep N_i$ is big, a contradiction. Thus, $P_{\CKM}\cdot
N_i=0$ for all $N_i$. Now suppose that $N^{\prime 2} > 0$, where
$N'\eqdef\sum\ep_iN_i$, $|\ep_i|\ll 1$. Then $(P_{\CKM}+N')
^2=N^{\prime 2} > 0$ and the divisor $P_{\CKM}+N'$ is big. Again
we have a contradiction.
\end{example}

The following theorem is a consequence of the existence of Iitaka
fibration and two-dimensional Zariski decomposition.

\begin{theorem}[{\cite{Cu}}]
Let $D$ be an effective Cartier divisor. Assume that $1\le\kappa
(X, D)\le 2$. Then there exists a birational contraction $f\colon
Y\to X$ such that $f^*D$ has an CKM-Zariski decomposition
$N_{\CKM} (f^*D)$ satisfying the following condition
\begin{equation}
\label{eq-kappa}
\kappa (X, D)=\nu (Y, P_{\CKM} (f^*D)).
\end{equation}
Furthermore, if $D$ is a $\QQ$-divisor, then so is $N_{\CKM} (D)$.
An CKM-decomposition satisfying condition \eqref{eq-kappa} is
unique.
\end{theorem}

It is expected that Zariski decomposition exists for the log
canonical divisor. In this instance Zariski decomposition must
have good properties:

\begin{theorem}[{\cite{Be}}, {\cite{Mo1}}, {\cite{Mo2}}, {\cite{Ka}}]
\label{th-CKM-logcanonical}
Let $(X,\Delta)$ be a projective Kawamata log terminal pair such
that $\Delta$ is a $\QQ$-boundary and the divisor $K_X+\Delta$ is
big. Assume that there exists a Zariski decomposition in CKM sense
for $D=K_X+\Delta$. Then the positive part $P_{\CKM} (D)$ is a
semiample $\QQ$-divisor. In particular, the algebra $\sR_X
(K_X+\Delta)$ is finitely generated.
\end{theorem}
An CKM-decomposition in this case coincides with Fujita and
Shokurov decompositions.

Theorem \xref{th-CKM-logcanonical} was proved by Moriwaki
{\cite{Mo1}}-{\cite{Mo2}} and Kawamata {\cite{Ka}} in more
general, relative situation. Thus the existence of a (relative)
Zariski decomposition in CKM sense for the log canonical divisor
of small contractions $X\to Z$ is a sufficient condition for the
existence of log flips.

\section{Decompositions of b-divisors}
\begin{pusto}
\label{def-pseudo-Mov-Fix}
Let $\sD$ be a b-divisor. Similar to \eqref{eq-Mov-Fix} we put
\begin{equation*}\Mov (\sD)\eqdef
\begin{cases}
-\inf\limits_{s\in K (X)^*}\{(s)\mid\sD+ (s)\ge 0\} &\text{if $H^0
(\Oh (\sD))\neq 0$},
\\
-\infty &\text{otherwise}.
\end{cases}
\end{equation*}
In the case $\Mov (\sD)\neq -\infty$, we may take a section
$s_0\in K (X)^*$ so that $\sD+ (s_0)\ge 0$. Then $-(s_0)\le\Mov
(\sD)\le D$. Therefore $\Mov (\sD)$ is a b-divisor. If $H^0 (\Oh
(\sD))\neq 0$, then we also put
\[
\Fix (\sD)\eqdef\inf\{\sL\mid\sL\sim\sD,\quad\sL\ge 0\}=\sD -\Mov
(\sD).
\]
B-divisors $\Mov (\sD)$ and $\Fix (\sD)$ defined above are called
\emph{mobile} and \emph{fixed} parts of $\sD$ respectively. If a
b-divisor $\sD$ is effective, then so is $\Mov (\sD)$: $\Mov
(\sD)\ge - (\Const) =0$. It is easy to see from the definition
that

\begin{equation}
\label{eq-Mov-Fix-b-div}
(\Mov (\sD))_X\le\Mov (\sD_X),\qquad (\Fix (\sD))_X\ge\Fix
(\sD_X).
\end{equation}
\end{pusto}

\begin{lemma}
\label{L-mobile-1-b}
Let $\sD$ be a b-divisor such that $H^0 (\Oh (\sD))\neq 0$. Then
$\sM=\Mov (\sD)$ satisfies the following properties:
\begin{enumerate}
\item
$\sM\le\sD$;
\item
the b-divisor $\sM$ is b-free \textup(in particular it is
b-Cartier\textup);
\item
if a b-divisor $\sL\le\sD$ is b-free, then $\sL\le\sM$.
\end{enumerate}
Conversely, if an \textup(integral\textup) b-divisor $\sM$
satisfies conditions {\rm (i) - (iii)}, then $\sM=\Mov (\sD)$.
\end{lemma}
\begin{proof}
Prove (ii). Similar to Lemma \xref{L-mobile-1} we have
\[
\Fix (\sM)=\inf\{\sL\mid\sL\sim\sM,\quad\sL\ge 0\}=0.
\]
Hence, $\Fix (\sM_X) =0$, i.e., the linear system $\lins{\sM_X}$
has no fixed components. Let $f\colon Y\to X$ be a resolution of
base points of the linear system $\lins{\sM_X}$ and let $\lins{M}$
be the birational transform of $\lins{\sM_X}$ on $Y$. Using
\PLF{Prop. 3.20} we obtain
$\lins{\sM_X}\simeq\lins{\sM_Y}\simeq\lins{M}$. Since
$\lins{\sM_Y}$ has no fixed components, $\lins{\sM_Y}=\lins{M}$
and this linear system has no base points. The rest of the proof
is completely similar to Lemma \xref{L-mobile-1}
\end{proof}

Similar to Lemma \xref{L-mobile-2} one can prove the following

\begin{lemma}
\label{L-mobile-2-b}
Let $\sD$ be a b-divisor such that $H^0 (\Oh (\sD))\neq 0$. Then
$\sM=\Mov (\sD)$ satisfies the following properties:
\begin{enumerate}
\item
$\sM\le\sD$;
\item
$H^0 (\Oh (\sM)) =H^0 (\Oh (\sD))$;
\item
if for a b-divisor $\sL\le\sD$ the equality $H^0 (\Oh (\sL)) =H^0
(\Oh (\sD))$ holds, then $\sL\ge\sM$.
\end{enumerate}
Conversely, if an \textup(integral\textup) b-divisor $\sM$
satisfies conditions {\rm (i) - (iii)}, then $\sM=\Mov (\sD)$.
\end{lemma}
Here (as well as everywhere) we assume that sections of the sheaf
$\Oh (\sD)$ are elements of the field $K (X)$:
\[
H^0 (\Oh (\sD))=\{s\in K (X)\mid\sD+ (s)\ge 0\}.
\]
Therefore the statement remains to be true for infinitely
dimensional spaces $H^0 (\Oh (\cdot))$.

\begin{definition}
A b-divisor $\sG$ is said to be \emph{pbs-ample} if there exists a
sequence b-semiample b-divisors $\sG_i$ such that
$\lim_{i\to\infty}\sG_i=\sG$.
\end{definition}

\begin{definition}
\label{prop-Zar_Shok-b}
A decomposition of a b-divisor $\sD=\sD^{\m}+\sD^{\e}$ is called a
\emph{pbs-decomposition} (pseudo b-semiample decomposition) if it
satisfies the following properties:
\begin{enumerate}
\item
$\sD^{\e}\ge 0$;
\item
$\sD^{\m}$ is pbs-ample;
\item
for any pbs-ample b-divisor $\sL$ such that $\sD^{\m}\le\sL\le\sD$
we have $\sL=\sD^{\m}$.
\end{enumerate}
\end{definition}

Here $\sD^{\m}$ is the \emph{maximal pbs-ample} part and
$\sD^{\e}$ is the \emph{fixed part}.

Pbs-decompositions generalize divisorial Shokurov decompositions~
\xref{def-Zar_Shok}. Indeed, let $D=D^{\m}+D^{\e}$ and $\sD^{\m}$
be a decomposition and the corresponding b-divisor from
\xref{def-Zar_Shok} and let $\sD^{\e}=D^{\e}$ (the last is
considered as an equality of b-divisors). Put also
$\sD=\sD^{\m}+\sD^{\e}$. Then by definition
$\sD=\sD^{\m}+\sD^{\e}$ is a pbs-decomposition and $D=\sD_X$.

\begin{pusto}
\label{def-s-decomp-b}
Similar to the divisorial case we say that a b-divisor $\sD$ is
\emph{effective modulo $\QQ$-linear equivalence} if $\sD\qq\sD'$,
where $\sD'\ge 0$. Equivalent: $H^0 (\Oh (n_0\sD))\neq 0$ for some
$n_0\in\NN$. Let $\sD$ be an effective modulo $\QQ$-linear
equivalence b-divisor. Denote
\[
\sM_n\eqdef\Mov (n\sD).
\]
Since $\sM_n/n\le\sD$, there exists the limit
\[
\sP_{\s} (\sD)\eqdef\limsup_{n\to\infty} (\sM_n/n).
\]
Since $\sM_{n_1n_2}\ge n_2\sM_{n_1}$ for all $n_1, n_2\in\NN$, we
have
\[
n\sP_{\s} (\sD)\ge\sM_n\quad\text{for all $n\in\NN$.}
\]

It is clear that $-(s_0)\le\sP_{\s} (\sD)\le\sD$, where $0\neq
s_0\in H^0 (\Oh (n_0\sD))$. Hence, $\sP_{\s} (\sD)$ is also a
b-divisor. Thus we obtain a decomposition
\[
\sD=\sP_{\s} (\sD)+\sN_{\s} (\sD),\qquad\sN_{\s} (\sD)\eqdef\sD
-\sP_{\s} (\sD)\ge 0,
\]
which we call an \emph{s-decomposition} of a b-divisor. Similar to
\eqref{eq-Ns-} we have
\begin{equation}
\label{eq-Ns-b-div}
\sN_{\s} (\sD)=\inf\{\sL\mid\sL\qq\sD,
\\
sL\ge 0\}.
\end{equation}
By construction $\sP_{\s} (\sD)$ is a pbs-ample b-divisor.
\end{pusto}

It is easy to see that properties \xref{svoistva} hold also for
the b-divisorial s-decomposition (if one replaces divisors with
b-divisors).

The concept of s-decomposition of b-divisors allows us to
formulate a criterion of finite generation of a \emph{b-divisorial
algebra} $\sR_X\sD\eqdef\oplus_{n\ge 0} H^0 (X,\OOO_X (n\sD))$
(cf. Theorem \xref{th-Lim-Crit}).

\begin{theorem}[Limiting Criterion \PLF{Th. 4.28}]
\label{th-Lim-Crit-b-div}
Let $\sD$ be an effective modulo $\QQ$-linear equivalence
b-divisor. Then the b-divisorial algebra $\sR_X\sD$ is finitely
generated if and only if $\sP_{\s} (\sD)=\Mov (n_0\sD)/ n_0$ for
some $n_0\in\NN$.
\end{theorem}

\begin{pusto}
Example \xref{ex-cremona} shows that the divisorial
s-decomposition $D=P_{\s} (D) +N_{\s} (D)$ does not agree via
pull-backs $f_*$. Therefore the system $P_{\s} (\sD_X)$ does not
form a b-divisor. In particular, $\sP_{\s} (\sD)_X\neq P_{\s}
(\sD_X)$. However, from \eqref{eq-Mov-Fix-b-div} we always have
\[
\sP_{\s} (\sD)_X\le P_{\s} (\sD_X),\qquad\sN_{\s} (\sD)_X\ge
N_{\s} (\sD_X).
\]
Where equalities are achieved ``birationally asymptotically'':
\end{pusto}

\begin{lemma}
\label{lemma-Ns=Ns}
\[
\sP_{\s} (\sD)_X=\inf_{f\colon Y\to X} f_*P_{\s}
(\sD_Y),\qquad\sN_{\s} (\sD)_X=\sup_{f\colon Y\to X} f_*N_{\s}
(\sD_Y),
\]
where the infimum and supremum are taken over all birational
contractions $f\colon Y\to X$.
\end{lemma}
\begin{proof}
Consider the following sets of divisors on $X$
\begin{eqnarray*}
&\mathfrak S\eqdef\{\sL_X\mid\sL\qq\sD,\quad\sL\ge 0\},
\\
&\mathfrak S_Y\eqdef\{\sL_X\mid\sL\qq\sD,\quad\sL_Y\ge 0\},
\end{eqnarray*}
where $Y$ is a birational model dominating $X$. Then $\mathfrak
S_{Y'}\subset\mathfrak S_Y$ if $Y'$ dominates $Y$ and $\mathfrak
S=\bigcap\limits_{Y/X}\mathfrak S_Y$. According to \eqref{eq-Ns-}
and \eqref{eq-Ns-b-div} we have $\sN_{\s} (\sD)_X=\inf\mathfrak S$
and $f_*N_{\s} (\sD_Y)=\inf\mathfrak S_Y$. Hence,
\[
\sN_{\s} (\sD)_X=\inf\mathfrak S=\inf\bigcap_{Y/X}\mathfrak
S_Y=\sup_{Y/X}\inf\mathfrak S_Y=\sup_{f\colon Y\to X} f_*N_{\s}
(\sD_Y).
\]
This proves the statement.
\end{proof}

It follows from definitions that $\sP_{\s} (\sD)\le\sD^{\m}$ for
any effective modulo $\QQ$-linear equivalence b-divisor $\sD$
(such that a decomposition $\sD=\sD^{\m}+\sD^{\e}$ exists).

Completely similar to \xref{remark-propert-Sh} one can prove the
following.
\begin{proposition}
\label{proposition-b-divisors}
Suppose that $\sP_{s} (\sD)$ is a rational b-semiample b-divisor.
Then for any b-semiample b-divisor $\sL\le\sD$ the inequality
$\sL\le\sP_{s} (\sD)$ holds.
\end{proposition}

According to Limiting Criterion \xref{th-Lim-Crit-b-div} the
conditions of proposition are satisfied, for example, in case when
the b-divisorial algebra $\sR_X\sD$ is finitely generated.

In conclusion we mention an interesting result of Fujita:

\begin{theorem}[{\cite{F5}}]
Let $D$ be a big Cartier divisor on a $d$-dimensional variety $X$
and let
\[
v (D)\eqdef\limsup_{t\to\infty}\left (\frac{d!}{t^n}\dim H^0
(X,\Oh (tD))\right)
\]
be its \emph{volume}. Then for any $\ep > 0$, there exists a
birational contraction $f\colon Y\to X$ and a decomposition
\[
f^*D= P_{\ep}+N_{\ep}
\]
in a sum of $\QQ$-Cartier divisors, where
\begin{itemize}
\item
$N_{\ep}\ge 0$;
\item
the divisor $P_{\ep}$ is semiample;
\item
$v (D) -\ep < v (P_{\ep})=(P_{\ep})^d < v (D)$.
\end{itemize}
\end{theorem}

\section{Analytic Zariski decomposition}
In conclusion we mention about complex analytic approach to the
constructing of Zariski decompositions \cite{Ts1}, \cite{Ts2} (see
also \cite{Bou}). Main definitions and facts on complex currents
can be found in \cite{GH} or \cite{De}.

Let $T$ be a closed positive $(1,1)$-current on the open unit ball
in $\CC^n$ with center at $0$. The \emph{Lelong number} $\Theta
(T, 0)$ of $T$ at $0$ is the number
\[
\Theta (T, 0)\eqdef\lim_{r\to 0}\frac{T (\chi (r)\omega^{n-1})}
{\pi^{n-1}r^{2 (n-1)}}
\]
where $\omega=\frac{\sqrt{-1}}{2}\sum dz_i\wedge d\bar z_i$ and
$\chi (r)$ is the characteristic function of the open ball in
$\CC^n$ with center at $0$ and radius $r$. This number is
invariant under the changes of coordinates. Thus it is possible to
define the Lelong number $\Theta (T, x)$ for a closed positive
$(1,1)$-current at a point $x\in X$ on any complex variety $X$. If
$T$ is given by integration over an analytic subvariety of
codimension $1$, then $\Theta (T, x)$ is the usual multiplicity of
this subvariety at $x$.

Further we suppose that $X$ is a nonsingular projective complex
variety. For any analytic subset $V\subset X$, we can define
\[
\Theta (T, V)\eqdef\inf\{\Theta (T, x)\mid x\in V\}.
\]
If $V$ irreducible, then $\Theta (T, V)$ coincides with the Lelong
number $\Theta (T, x)$ at a \emph{very} general point $x\in V$.
The last enable us to define \emph{Siu decomposition} (see
\cite[Ch. III, (8.16)]{De}):
\begin{equation}
\label{eq-AZD-Siu}
T=T'+\sum_{V}\Theta (T, V) [V],
\end{equation}
where the (infinite in general) sum is taken over all prime
divisors $V$ and $[V]$ is current given by integration over $V$.
Here $T'$ is also a closed positive $(1,1)$-current such that the
set $\{x\in X\mid\Theta (T, x)\ge\ep\}$ has codimension $\ge 2$
for any $\ep> 0$.

Let $f\colon Y\to X$ be a surjective morphism and let $T$ be a
closed positive $(1,1)$-current on $X$. Locally we may write
$T=\sqrt{-1}\partial\bar\partial\var$, where $\var$ is a
plurisubharmonic function \cite[Ch. III, Prop. 1.19] {De}.
Therefore we can define the pull-back of a current $T$ by the
formula $f^*T=\sqrt{-1}\partial\bar\partial f^*\var$ (this
definition does not depend on the choice of the function $\var$).

It is possible also to define a ``b-divisorial'' version of the
Siu decomposition: for a closed positive $(1,1)$-current $T$ we
put
\begin{equation}
\label{eq-AZD-Siu-1}
\sS (T)=\sum_{V}\Theta (T, V) V,
\end{equation}
where the sum is taken over all divisorial discrete valuations $V$
of the field $K (X)$ and $\Theta (T, V)$ is computed on a
(nonsingular projective) birational model $f\colon Y\to X$, on
which the center of $V$ is a prime divisor: $\Theta (T, V)=\Theta
(f^*T, V)$. Strictly speaking $\sS (T)$ is not a b-divisor, since
$\sS (T)_X$ can be an infinite sum. Nevertheless, $\sS (T)$
satisfies the condition $\sS (T)_{Y_1}=g_*\sS (T)_{Y_2}$ whenever
there is a birational morphism $g\colon Y_1\to Y_2$.

\begin{theorem}[\cite{Ts1}]
Let $D$ be a big divisor on $X$. Then there exists a
\textup(closed, positive\textup) $(1,1)$-current $T$ such that
\begin{enumerate}
\item
$T$ represents the cohomology class $c_1 (D)$ of $D$ in $H^2
(X,\RR)$;
\item
for any birational contraction $f\colon Y\to X$, any nonnegative
integer $m$, and any point $y\in Y$ the following approximation
formulas for $\Theta (f^*T, y)$ hold:
\begin{eqnarray}
\label{eq-AZD-1}
&&\min_{L\in |f^*mD|}\left (\mult_y L\right)\ge m\Theta (f^*T, y)
\\
\label{eq-AZD-2}
&&\liminf_{m\to\infty}\left (\frac1m\min_{L\in |f^*mD|} \left
(\mult_y L\right)\right)=\Theta (f^*T, y).
\end{eqnarray}
\end{enumerate}
\end{theorem}

Such a current $T$ is called an \emph{analytic Zariski
decomposition \textup(AZD\textup)} of $D$.

Equality \eqref{eq-AZD-2}, in particular, asserts that the Siu
decomposition of a current $T$ is a divisor and coincides with the
s-decomposition $N_{\s} (D)$. Furthermore, the formal sum $\sS
(T)$ is a b-divisor and $\sS (T)=\sN_{\s} (\ov D)$.

\begin{corollary}
Let $D$ be a divisor on $X$.
\begin{enumerate}
\item
If $D$ is nef and big, then $c_1 (D)$ can be represented by a
closed positive $(1,1)$-current $T$ on $X$ such that $\Theta
(T)\equiv 0$.
\item
Conversely, if $c_1 (D)$ can be represented by a closed positive
$(1,1)$-current $T$ on $X$ such that $\Theta (T)\equiv 0$, then
$D$ is nef.
\end{enumerate}
\end{corollary}

\par\bigskip\noindent
The author would like to thank Professors V.\, V. ~Shokurov and
V.\, A. ~Iskovskikh for useful discussions as well as for that
they looked through the preliminary manuscript and gave numerous
remarks. A part of this work was carried out at
Max-Planck-Institut f\"ur Mathematik in 2001 and at the Isaac
Newton Institute for Mathematical Sciences in 2002. Author very is
grateful to the research staff and administrations of these
institutes for hospitality and wonderful working environment.

\end{document}